\documentclass[12pt]{article}
\usepackage{latexsym}
\usepackage{amsmath}
\usepackage{amssymb}
\usepackage{wasysym}
\usepackage{bbm}
\usepackage{color}
\usepackage{graphicx}
\usepackage{epsfig}
\usepackage[tmargin=1.25in,bmargin=1.25in,lmargin=1.1in,rmargin=1.1in]{geometry}
\newtheorem{thm}{Theorem}[section]
\newtheorem{lem}[thm]{Lemma}

\newtheorem{cor}[thm]{Corollary}

\newtheorem{rem}[thm]{Remark}
\makeatletter\@addtoreset{equation}{section}\makeatother
\newcommand{\ea}{\end{array}}
\newcommand{\beqohne}{\begin{eqnarray*}}
\newcommand{\eeqohne}{\end{eqnarray*}}
\newcommand{\beohne}{\begin{equation*}}
\newcommand{\eeohne}{\end{equation*}}
\newcommand{\R}{\mathbb{R}}

\newcommand{\C}{\mathbb{C}}
\def\proof{\noindent{\bf Proof:}\hskip10pt}

\def \sur#1#2{\mathrel{\mathop{\kern 0pt#1}\limits^{#2}}}

\newcommand{\mun}{\mu^{(n)}}
\newcommand{\muun}{\mu^{(n)}_{\u}}

\def \u{{\tt u}}

\def \v{{\tt v}}
\def \ben{\begin{eqnarray}}
\def \een{\end{eqnarray}}
\newcommand{\tr}{\mathrm{tr}}
\newcommand{\Dir}{\operatorname{Dir}}
\newcommand{\ext}{{\operatorname{ext}}}

\newcommand{\SC}{\operatorname{SC}}
\newcommand{\Pnv}{\mathbb P_n^V}

\newcommand{\muunI}{\mu^{(n)}_{\u,I}}

\newcommand{\sn}{^{(n)}}
\newcommand{\Ir}{\mathcal{J}}
\newcommand{\Fr}{\mathcal{F}}
\newcommand{\Sr}{\mathcal{S}}

\newcommand{\diff}{\mathop{}\!\mathrm{d}}

\definecolor{Red}{rgb}{1,0,0}
\definecolor{Blue}{rgb}{0,0,1}
\def\red{\textcolor{Red}}

\begin{document}

\title{Sum rules via large deviations: extension to polynomial potentials and the multi-cut regime}

\author{{\small Fabrice Gamboa}\footnote{ Universit\'e Paul Sabatier, Institut de Math\'ematiques de Toulouse,  31062-Toulouse Cedex 9, France and ANITI, 
{\tt fabrice.gamboa@math.univ-toulouse.fr}}
\and{\small Jan Nagel}\footnote{Technische Universit\"at Dortmund, Fakult\"at f\"ur Mathematik, 44227 Dortmund, Germany,
{\tt jan.nagel@tu-dortmund.de}}
\and{\small Alain Rouault}\footnote
{Laboratoire de Math\'ematiques de Versailles, UVSQ, CNRS, Universit\'e Paris-Saclay, 78035-Versailles Cedex France, {\tt alain.rouault@uvsq.fr}}}

\maketitle

\begin{abstract}
A sum rule is an identity connecting the entropy of a measure with coefficients involved in the construction of its orthogonal polynomials (Jacobi coefficients).  
Our paper is an extension of \cite{magicrules}
where we have showed sum rules by using only probabilistic tools (namely the large deviations theory). 
Here, %In this new paper, 
we  prove  large deviation principles for the weighted spectral measure of unitarily invariant random matrices in two general situations: firstly, %. First 
when the equilibrium measure is not necessarily supported by a single interval and secondly, when the potential is a nonnegative polynomial. %potentials expressing 
The rate functions can be expressed as  functions of the Jacobi coefficients.
 These new large deviation results lead to original sum rules both for the one and the multi-cut regime and also answer a conjecture stated in \cite{magicrules} concerning general sum rules.  
\end{abstract}

%%%%%%%%%%%%%%%%%%%%%%%%%%%%%%%%%%%%%%%%%%%%%%%%%%%%%%%%%%%%%%%%%%%%%%%%%
%%%%%%%%%%%%%%%%%%%%%%%%%%%%%%%%%%%%%%%%%%%%%%%%%%%%%%%%%%%%%%%%%%%%%%%%%


%\maketitle

%%%%%%%%%%%%%%%%%%%%%%%%%%%%%%%%%%%%%%%%%%%%%%%%%%%%%%%%%%%%%%%%%%%%%%%%%
%%%%%%%%%%%%%%%%%%%%%%%%%%%%%%%%%%%%%%%%%%%%%%%%%%%%%%%%%%%%%%%%%%%%%%%%%

\section{Introduction}
\label{sec:introduction}

This paper deals with the so-called sum rules arising from  spectral theory and orthogonal polynomials on the real line (OPRL). Our approach uses only probabilistic methods. 
Given a probability measure $\mu$  with compact support on $\R$, we may encode %parametrize
 $\mu$ by the recursion coefficients of orthonormal polynomials in $L^2(\mu)$. A sum rule is an identity between a non-negative functional of these coefficients and an entropy-like functional of $\mu$, 
each side giving the discrepancy between $\mu$ and some reference measure. 
 
When the reference measure is the semicircle distribution
\begin{equation}
\label{SC0}
\operatorname{SC}(dx) = \frac{1}{2\pi}\sqrt{4-x^2}\!\ \mathbbm{1}_{[-2, 2]}(x)\!\ dx, 
\end{equation}  
the sum rule was proved with spectral theory method by Killip and Simon in \cite{KS03}\footnote{An exhaustive discussion and history of this sum rule can be found in  Section 1.10 of the book  \cite{simon2} and a deep analytical proof is in Chapter 3.
}. This result is the OPRL counterpart of the classical Szeg{\H o} theorem for orthogonal polynomials on the unit circle (OPUC), where the reference measure is the Lebesgue measure.
An important consequence of such equalities is the equivalence of two conditions for the finiteness of both sides,  one formulated in terms of Jacobi coefficients and the other as a spectral condition. In the words of Simon \cite{simon2}, these are the \emph{gems} of spectral theory.

In \cite{GNROPUC}  and \cite{magicrules}, we 
gave a probabilistic interpretation of these sum rules and a general strategy to construct and prove new sum rules. 
In the OPRL case, the starting point is a random $n\times n$ Hermitian matrix $X_n$ and a fixed vector $e \in \C^n$. The random spectral measure $\mu_n$ of the pair $(X_n , e)$ is a weighted sum of Dirac masses supported by the (real) eigenvalues of $X_n$.
When the density of $X_n$ is unitarily invariant,  proportional to $\exp\{ -\frac{\beta}{2} n  \tr\!\ V(X)\}$ with a confining potential $V$,  we proved that $\mu_n$ satisfies the 
Large Deviation Principle (LDP) with speed $n$ and good rate function $\mathcal I_{\mathrm{sp}}$ involving the reversed entropy with respect to a measure $\mu_V$. This equilibrium measure $\mu_V$, or reference measure is the minimizer of the rate function or equivalently, the limit of the spectral measure as $n\to\infty$.  
Besides, in all the classical ensembles (Gaussian, Laguerre and Jacobi ensemble), the random recursion coefficients have a nice probabilistic structure (independence or slight dependence) so that we proved also an LDP for the ``coefficient encoding" of $\mu_n$ with speed $n$ and rate function $\mathcal I_{\mathrm{co}}$. Since a large deviation rate function is unique,
this implies the identity $\mathcal I_{\mathrm{sp}} = \mathcal I_{\mathrm{co}}$. For the Gaussian ensemble, this identity is precisely the sum rule of Killip and Simon. 
%The randomization with
For  the Laguerre or Jacobi ensemble it leads to new sum rules, with %respectively
reference measures the Marchenko-Pastur and the Kesten-McKay distributions, respectively. Furthermore, this method could be generalized to measures on the unit circle \cite{GNROPUC} or to operator valued measures \cite{GaNaRomat,GaNaRoJac}. Besides, it provides evidence for the Lukic conjecture \cite{BSZ}, see also \cite{BSZ1} for an exposition of the method.

For the measure side, the common feature of these models is the assumption that 
the equilibrium measure is supported by a single compact interval. In statistical physics terms this is the one-cut case, in contrast to the multi-cut case when the support is a finite union of disjoint compact intervals. In spectral theory, the first situation is called ``no gap" and the second one ``a finite number of open gaps".

For the coefficient side, the common feature is sufficient stochastic independence of the Jacobi coefficients. Nevertheless in Section 3.3 of \cite{magicrules},  based on \cite[Proposition 2]{KriRidVir2016}, we conjectured that, under some suitable conditions on $V$, the rate function on the coefficient side could be an expression with 
some limit involving $\tr\!\ V(T_n)$ as $n\to \infty$, where $T_n$ is the $n$-dimensional Jacobi matrix. 

Besides, by spectral theory methods, \cite{Nazarov} obtained a more general sum rule, when the reference measure is $A(x) \operatorname{SC}(dx)$ with $A$ a nonnegative polynomial (see the discussion in Section \ref{suonecut}). This is equivalent to start from a one-cut polynomial $V$.

Here, we extend our probabilistic method along two directions. 
Firstly, we  show a large deviation theorem for the spectral measure sequence $(\mu_n)_n$ in the multi-cut case, for general potentials $V$ %with equilibrium measure supported by a finite number of compact intervals
 (Theorem \ref{MAIN}). Secondly, when $V$ is a nonnegative polynomial 
we show an  LDP in terms of the Jacobi coefficients (Theorem \ref{MAIN2}).
Surprisingly,  the rate function in this new LDP  contains a %new
 remainder term, which actually vanishes in the  case of a polynomial potential with one-cut equilibrium measure.  These  two last results are obtained by a similar method as developed by Breuer, Simon and Zeitouni \cite{BSZ} (for a polynomial potentials in the OPUC case with full support equilibrium measure). Indeed, the crucial argument  to tackle the remainder term is the Rakhmanov's theorem (see \cite{rahmanov1977asymptotics}, \cite{Den2004}). 

The combination of our new LDPs leads to a general gem in the multi-cut polynomial case, (Theorem \ref{abstractgem}), and an exact sum rule in the one-cut-polynomial case, (Theorem \ref{newsumrule}).   
While convex potentials lead to a one-cut equilibrium measures, the new gem also applies to nonconvex polynomial potentials.   
We guess that the new sum rule may hold true for more general potentials including in particular one or two logarithmic contribution(s). 
In \cite{magicrules} Sec. 3.3.1, it is proved that for  Laguerre and Jacobi potentials the claim of Theorem  \ref{newsumrule} holds true.

Other gems, i.e., sets of equivalent conditions for spectral measures in the multi-cut case were given by \cite{eichinger2016jacobi} and \cite{yuditskii2018killip} based on the Jacobi flow approach. Our method yields another expression for the coefficient side which depends on the potential in a natural way. The different formulations illustrate the different point of views: whereas the spectral theoretic methods start from a perturbation of the semicircle law (or free Jacobi matrix), the starting point for our probabilistic approach is a randomization given by the potential $V$.

\section{Notations and definitions}
\label{sec:sumrules}

\subsection{Tridiagonal representations}

Let $\mathcal{M}_1$ be the set of all probability measures on $\mathbb{R}$.  
%Given a probability measure
For $\mu\in\mathcal{M}_1$ with compact but infinite support (known as the nontrivial case), let $p_0,p_1,\dots $ be the orthonormal polynomials with positive leading coefficients obtained by applying the orthonormalizing Gram-Schmidt procedure
to the sequence $1, x, x^2, \dots$ in $L^2(\mu)$. They obey the recursion relation
\begin{align} \label{polrecursion}
xp_k(x) = a_{k+1} p_{k+1}(x) + b_{k+1} p_k (x) + a_{k} p_{k-1}(x)
\end{align}
for $ k \geq 0$ (resp. for $0 \leq k \leq n-1$) where the Jacobi parameters satisfy $b_k \in \mathbb R, a_k > 0$ for all $k\geq 1$ and with $p_{-1}(x)=0$.

 In the basis $\{p_0, p_1, \dots\}$, the 
linear transform $f(x) \rightarrow xf(x)$ (multiplication by the identity) 
 in $L^2(d\mu)$ is represented by the matrix
\ben
\label{favardinfini}
T = \begin{pmatrix} b_1&a_1 &0&0&\cdots\\
a_1&b_2  &a_2&0&\cdots\\
0&a_2 &b_3&a_3& \\
\vdots& &\ddots&\ddots&\ddots
\end{pmatrix} ,
\een
where we have $a_k > 0$ for every $k$. The mapping $\mu \mapsto T$ (called here the Jacobi mapping) is a one to one correspondence between probability measures on $\mathbb R$ having compact infinite support and Jacobi matrices 
with built with sequences satisfying $\sup_n(|a_n| +
|b_n|) < \infty$. Actually, such Jacobi matrix is identified as an element of $\mathcal{R}$ (defined below in \ref{defR}). This result is 
sometimes called Favard's theorem. % (see \cite{Simon3} p. 432). 
If $T$ is a infinite Jacobi matrix, we denote the $N\times N$ upper left subblock  by $\pi_N(T)$.

If $\mu\in \mathcal{M}_1$ is supported by $n$ distinct points, we may still define orthonormal polynomials up to degree $2n-1$. The $n$-th polynomial has roots at the $n$ support points and thus norm zero in $L^2(\mu)$. We then consider the finite dimensional Jacobi matrix of $\mu$, 
\ben
\label{favardfini}
T_n = \begin{pmatrix} b_1&a_1 &0&\dots&0\\
a_1&b_2  &a_2&\ddots&\vdots\\
0&\ddots &\ddots&\ddots&0\\
\vdots&\ddots&a_{n-2}&b_{n-1}&a_{n-1}\\
0&\dots&0&a_{n-1}&b_{n}
\end{pmatrix} .
\een
So, measures supported by $n$ points lead to $n\times n$ symmetric tridiagonal matrices with subdiagonal positive terms. In fact, there
is a one-to-one correspondence between such a matrix and such a measure. 
We can identify $T_n$ 
with the vector $r_n=(b_1,a_n\dots ,a_{n-1},b_n)$. It is convenient to embed this into sequence spaces and to
 identify $r_n$ with the infinite vector  $(b_1,a_1\dots ,a_{n-1},b_n,0,\dots)\in \mathcal{R}$, where 
\begin{align}
\label{defR}
\mathcal{R} = (\mathbb{R}\times [0,\infty))^{\mathbb{N}}\,.
\end{align} 
Similarly, $T_n$ may be identified with the one-sided infinite Jacobi matrix extended by zeros. %In this setting, we regard $r_n$ as a random element of  $R$, equipped with the product topology and the Borel $\sigma$-algebra.
For an element $r\in \mathcal{R}$, we let $\pi_N(r)\in \mathcal{R}_N$ be the projection onto the first $2N-1$ coordinates.

Let now $\psi$ be
 the mapping defined on the set of measures $\mu\in \mathcal{M}_1$ with compact support by
\begin{align}
\psi(\mu) = r = (b_1,a_1,b_2,\dots ) ,
\end{align}
 where $a_n=b_{n+1}=0$ if $\#\operatorname{supp}(\mu)\leq n$. Note that $\psi$ is not continuous. Nevertheless, if for $K > 0$ 
\begin{align}
\label{defMK}
\mathcal{M}_{1,K} = \{\mu \in \mathcal M_1 : \operatorname{supp} (\mu) \in [-K, K] \}\,,
\end{align}
then, $\psi$ is a homeomorphism on $\mathcal{M}_{1,K}$.

An other point of view consists in considering that the measure $\mu$ is the spectral measure of the tridiagonal operator. 
More precisely, let $H$ be a self-adjoint bounded operator on a Hilbert space $\mathcal H$ and $e\in \mathcal H$ be a cyclic vector (that is, such that the linear combinations of the sequence $(H^ke)$ are dense in $\mathcal H$). 
Then, the spectral measure of the  pair $(H,e)$ is the unique 
$\mu\in \mathcal M_1$
 such that
\[\langle e, H^k e\rangle = \int_\mathbb R x^k \diff \mu(x) \ \ (k \geq 1).\]
It turns out that $\mu$ is a unitary invariant for $(H,e)$. Another invariant is the tridiagonal reduction 
defined above.

 If dim $\mathcal H =n$ and $e$ is cyclic for $H$,  let $\lambda_1, \ldots, \lambda_n$ be the (real) eigenvalues of $H$ and let $u_1, \ldots, u_n$ be a system of orthonormal eigenvectors. The spectral measure of the pair $(H,e)$ is then
\begin{align}\label{spectralmeasure}
\mun =  \sum_{k=1}^n w_k\delta_{\lambda_k}\,,
\end{align}
with $w_k= |\langle u_k, e\rangle|^2$. 
This measure is a weighted version of the empirical eigenvalue distribution
\begin{align}\label{empiricallaw}
\muun = \frac{1}{n} \sum_{k=1}^n \delta_{\lambda_k} \,.
\end{align}

If $J$ is  a Jacobi matrix, 
we can take the first vector $e_1$ of the canonical basis as the cyclic vector. 
 Let  $\mu$  be the spectral measure associated to the pair $(J, e_1)$, then $J$ represents the multiplication by $x$ in the basis of orthonormal polynomials associated to $\mu$ and $J=T(\mu)$.

Although  the general sum rules of the present work are purely deterministic identities, we need now to 
present a randomization to define  the elements involved in these formulas.

\subsection{Randomization}

\label{sec:randomization}
In the following $\beta = 2\beta' >0$ is a parameter, having in statistical physics the meaning of inverse temperature. %and we set $\beta' = \frac{\beta}{2}$.

The main object in our large deviation results is the random probability measure 
\begin{align} \label{defmun}
\mu_n = \sum_{k=1}^n w_k\delta_{\lambda_k} .
\end{align}
For suitable $V$, we let $\mathbb{P}_n^V$ be the distribution of a random measure
 such that
\begin{itemize}
\item the support points 
$(\lambda_1,\dots ,\lambda_n)$ have the joint density
\begin{align}\label{generaldensity}
 (Z_n^V)^{-1} 
e^{- n\beta' \sum_{k=1}^nV(\lambda_k)}\prod_{1\leq  i < j\leq n} |\lambda_i - \lambda_j|^\beta.
\end{align}
with respect to the Lebesgue measure on $\mathbb{R}^n$, 
\item
the weights $(w_1,\dots w_n)$ have a Dirichlet distribution $\Dir_n(\beta')$ of homogeneous parameter $\beta'$ 
 on the simplex $\{(w_1,\dots ,w_n)\in [0,1]^n |\, \sum_k w_k = 1\}$,  
with density proportional to $(w_1 \cdots w_n)^{\beta' - 1}$.
%which is defined by the Lebesgue density for the first $n-1$ coordinates proportional to 
%\begin{align*}
%\big( x_1 \cdots x_{n-1}(1-x_1-\dots -x_{n-1}) \big)^{\beta' -1} \mathbbm{1}_{\{ x_i > 0, x_1+\dots +x_{n-1}<1 \} } , 
%\end{align*}
%and additionally, 
\item the support points $(\lambda_1,\dots ,\lambda_n)$ are independent of the weights $(w_1,\dots w_n)$. 
\end{itemize}
Formula (\ref{generaldensity}) defines a log-gas density  of particles in an external potential $V$.
 
For specific values of $\beta$, %$\beta= 1, 2,4$,
 the distribution of $\mu_n$ is exactly the distribution of the spectral measure as defined in the above section.
For $\beta= 1$ (resp. $\beta= 2, \beta=4$), it is  the distribution of the spectral measure of the pair $(X_n , e_1)$ where $X_n$  is a random symmetric (resp. Hermitian, self-dual) matrix  whose density  is proportional to $\exp-\beta' V(X)$.

Additionally, for some classical potentials (Hermite, Laguerre, Jacobi) and general $\beta$, there is models of tridiagonal random matrices whose spectral measures are distributed as $\mu_n$ (see \cite{dumede2002}, \cite{Killip1}).
For general potentials and general $\beta$, it is shown in \cite{KriRidVir2016}, Proposition 2, that under $\Pnv$, the Jacobi coefficients $r_n=(b_1,a_1,\dots , b_n)$ have a density proportional to %density
\begin{align} \label{krishnadensity}
 (\tilde Z_n^V)^{-1} \exp \left\{ -n\beta' \left( \tr\!\ V(T_n) - 2\sum_{k=1}^{n-1} (1-\tfrac{k}{n} - \tfrac{1}{n\beta}) \log (a_k) \right) \right\} 
\end{align}
with respect to the Lebesgue measure on $\mathcal{R}_n = (\mathbb{R}\times [0,\infty))^{n-1}\times \mathbb{R}$ and where $T_n$ is as in \eqref{favardfini}. 

\iffalse
Notice that for very special values of $\beta$, the underlying probability distributions on matrices could be explicitly computed.
Indeed,  the corresponding random matrix is $n\times n$ self-adjoint real ($\beta=1$), complex ($\beta=2$) or quaternion ($\beta=4$). Furthermore, it has a
density proportional to $\exp -n\beta' V(X)$. Now, the spectral measure $\mu_n$ of a random matrix $X_n$ (at the vector $e_1$), is defined by the relation 
\begin{align} \label{spectralmeasuremoments}
\int x^k \diff \mu_n(x) = (X_n^k)_{1,1} ,\qquad k\geq 0,
\end{align}
and is given  in \eqref{defmun}.
(see for example \cite{mehta04}). 
The support points in this case are then precisely the eigenvalues of the matrix $X_n$. 
\red{So that, the similarity of $\mu^{(n)}$ in \eqref{spectralmeasure} and $\mu_n$ in \eqref{defmun} is not a coincidence!}
 is shown in \cite{KriRidVir2016}, Proposition 2, that under $\Pnv$, the Jacobi coefficients $r_n=(b_1,a_1,\dots , b_n)$ have a density proportional to %density
\begin{align} \label{krishnadensity}
 (\tilde Z_n^V)^{-1} \exp \left\{ -n\beta' \left( \tr\!\ V(T_n) - 2\sum_{k=1}^{n-1} (1-\tfrac{k}{n} - \tfrac{1}{n\beta}) \log (a_k) \right) \right\} 
\end{align}
with respect to the Lebesgue measure on $\mathcal{R}_n = (\mathbb{R}\times [0,\infty))^{n-1}\times \mathbb{R}$ \red{and where $T_n$ is as in \eqref{favardfini}.} 
\fi

\subsection{Assumptions on the potential}

The potential $V:\mathbb{R} \to (-\infty,+\infty]$ is supposed to be continuous and real valued on the interval ${(b^-, b^+)}$ ($-\infty\leq b^-<b^+\leq+\infty$), infinite outside of $[b^-,b^+]$ and $\lim_{x\to b^\pm} V(x) = V(b^\pm)$ with possible limit $V(b^\pm)=+\infty$. 
%Let $\beta=2\beta'>0$ be the inverse temperature. 
We will always make the following assumption. 
\begin{itemize}
\item[(A1)] Confinement: If $|b|=\infty$, $b\in\{b^-,b^+\}$, then
\begin{align*}
\liminf_{x \rightarrow b} \frac{V(x)}{2 \log |x|} > 1 .%\max(1, \beta^{-1}) . 
\end{align*}
\end{itemize}
Under (A1), the functional $\mathcal E (\mu)$ defined by
%the empirical distribution $\muun$ of eigenvalues $\lambda_1,\dots ,\lambda_n$ has a limit $\mu_V$ (in probability), which is compactly supported and the unique minimizer of %\eqref{ratemuu}
\begin{align}
\label{ratemuu}
\mu  \mapsto \mathcal E (\mu) := \int V(x) d\mu(x) - \int\!\!\!\int\log |x-y| d\mu(x)d\mu(y)
\end{align}
has a unique minimizer $\mu_V$ which is compactly supported, see \cite{johansson1998fluctuations} or \cite{agz}.  
We write $b_k^V,a_k^V$ for the Jacobi coefficients of $\mu_V$.
Further, we denote by $I$ the support of $\mu_V$.

The following assumption is crucial for the large deviation behavior of the extremal eigenvalues. 

\noindent (A2)  Control (of large deviations): the effective potential
\begin{align}
\label{poteff}
\Ir_V (x) := V(x) -2\int \log |x-\xi|\!\ d\mu_V(\xi)
\end{align}
achieves its global minimum value on $(b^-, b^+) \setminus \operatorname{Int}(I)$  only on the boundary of this set.

We need also the function 
\begin{align} \label{defF} 
\mathcal{F}_V(x) & = \begin{cases}
\mathcal{J}_V(x) - \inf_{\xi \in \R} \mathcal{J}_V(\xi) & \text{ if } x \notin \operatorname{Int}(I), \\
\infty & \text{ otherwise. } 
\end{cases}
\end{align}

For $d \in \mathbb N$, let $\mathcal V_{2d}$ the set of all polynomials of degree $2d$ with coefficient of the leading term positive, and 
\[\mathcal V = \bigcup_{d\geq 1} \mathcal V_{2d}\,.\]

It is known (\cite{pastur_shcherbina} Sect. 11.2) that if $V \in \mathcal V_{2d}$, then the support of $\mu_V$ is the union of a finite number $M \leq d$ of disjoint intervals
\[I = I_1\cup \dots \cup I_M\,.\]
If $M=1$, we say that we are in the one-cut regime, otherwise, we are in the multi-cut regime.
 
 Notice  that when $V\in \mathcal V$ is convex, then we are in the one-cut case and assumption (A2) is satisfied \cite[Proposition 3.1]{johansson1998fluctuations}.

\section{Sum rules}

Let $I = I_1\cup \dots \cup I_M$ be a union of $M$  compact and disjoint intervals, each with nonempty interior. We introduce the set $\Sr =\Sr(I)$ of finite non-negative measures $\mu$ which have compact support 
\begin{align} \label{support}
\operatorname{supp}(\mu) = J \cup E ,
\end{align}
with $J\subset I$ and $E=E(\mu)$ a finite or countable subset of $\mathbb{R}\setminus I$. 
% whose only points of accumulation are in $\partial I$. 
Denote by $\Sr_1(I)$ the set of all probability measures
member of $\Sr$. 
For a point $\lambda \in \mathbb{R}$, we denote by $e(\lambda)$ the point in $\partial I$ minimizing the distance to $\lambda$, where in case of ties we choose the leftmost one. 

In all the sum rules considered, the Kullback-Leibler divergence or relative entropy between two probability measures  
$\mu$ and $\nu$ plays a major role. When the ambient space is $\mathbb R$ endowed with its Borel  $\sigma$-field 
%$\mathbb{R}$),  %this  For $\mu$ and $\nu$ probability measures on $\mathbb R$, 
it is defined by
\begin{equation}
\label{KL}
{\mathcal K}(\mu\, |\ \nu)= \begin{cases}  \ \displaystyle\int_{\mathbb R}\log\frac{d\mu}{d\nu}\!\ d\mu\;\;& \mbox{if}\ \mu\ \hbox{is absolutely continuous with respect to}\ \nu ,\\
%\;\mu\ll \nu\;\mbox{ and }\log\frac{dP}{dQ}\in L^1(P),\\
   \  \infty  &  \mbox{otherwise.}
\end{cases}
\end{equation}
Usually, $\nu$ is the reference measure. Here the spectral side will involve the reversed Kullback-Leibler divergence, where $\mu$ is the reference measure and $\nu$ is the argument. In this case, we have that $\mathcal{K}(\mu |\nu)$ is finite if and only if
\begin{align} \label{eq:KL2}
\int \log w(x) \, d\mu(x) > - \infty , 
\end{align} 
where $d\nu(x) = w(x)d\mu(x) + d\nu_s(x)$ is the Lebesgue decomposition of $\nu$ with respect to $\mu$.

\subsection{A general gem}

The following theorem is our first spectral theoretic main result. It is a \emph{gem}, as explained in the introduction, giving equivalent conditions when a measure $\mu$ is sufficiently ``close'' to the reference measure $\mu_V$, defined through its potential $V$. We denote by $r^V=(b_1^V,a_1^V,b_2^V,\dots )$ the Jacobi coefficients of $\mu_V$.   

\begin{thm} \label{abstractgem}
Let $V\in \mathcal V$ %be a polynomial of even degree
 and let 
 $\mu$ be a probability measure with compact, infinite support. Then
\begin{align}
\label{gengem}
\sup_{N\geq 1} \left[ \tr V(\pi_N(r)) -\tr V(\pi_N(r^V)) - 2\sum_{k=1}^{N-1} \log (a_k/ a_k^V) \right] <\infty
\end{align} 
if and only if
\begin{enumerate}
\item 
  $\mu \in \mathcal S_1(I)$ with $I=\operatorname{supp}(\mu_V)$,
\item  $\sum_{\lambda \in E(\mu)} \mathcal{F}_V(\lambda)<\infty$ ,
%$\sum_{\lambda \in E}\, (\lambda - e(\lambda))^{3/2}  < \infty$  
\item the Lebesgue decomposition $\diff \mu(x) = f(x) \diff\mu_V(x)+\diff\mu_s(x)$ 
with respect to $\mu_V$
satisfies
\begin{align*}
\int_I \log f(x)\, \diff\mu_V(x) >-\infty .
\end{align*}
\end{enumerate}
\end{thm}

\begin{rem} \label{rem:summability}
\begin{enumerate}
\item In the vocabulary of spectral theory, Condition 1 is called Blumenthal-Weyl and Condition 3 is called Quasi-Szeg{\H o} condition.
\item
The second condition in Theorem \ref{abstractgem} regarding the outliers in $E(\mu)$ may be simplified, provided one can specify the decay of the density of $\mu_V$ at the boundary of $I$. Suppose the equilibrium measure $\mu_V$ has a Lebesgue-density $\rho_V$ satisfying $\rho_V(x) = \sqrt{d(x,\partial I)}Q(x)$ with $Q$ and $Q^{-1}$ bounded on $\operatorname{supp}(\mu_V)$ (see, e.g., Theorem 11.2.1 in \cite{pastur_shcherbina} for sufficient conditions), then condition $2.$ in Theorem \ref{abstractgem} is equivalent to   
\begin{align}
\label{Lieb}
\sum_{\lambda \in E}\, d(\lambda,I)^{3/2}  < \infty .
\end{align}
In the vocabulary of spectral theory, it is then called Lieb-Thirring condition.
\end{enumerate}
\end{rem}

\subsection{The one-cut sum rule}
\label{suonecut}

We now give a general sum rule for any suitable polynomial potential.
%Indeed, the following result shows that for any suitable polynomial potential a sum rule is available.
% Surprisingly, 
In this quite general frame the shape of the sum rule
remains  the same as previously. On one hand, the spectral side  always involves both the reversed Kullback information with respect to the corresponding equilibrium  measure  and an additional term related to the non-essential spectra. On the other hand, the other side is a discrepancy between the Jacobi coefficients. Hence, this Theorem gives rise to a 
general sum rule for the one-cut case. We will see later (in Theorem \ref{MAIN2}) that, in the multi-cut case a remainder term in the Jacobi coefficients side appears.  
%In the one-cut case%when $\mu_V$ is supported by a single interval
%, we can show the following sum rule.

\begin{thm} \label{newsumrule}
Let $V\in \mathcal V$ such that $\mu_V$ is supported by a single interval and let $\mu$ be a probability measure with compact  infinite support. Then
\begin{align*}
\lim_{N\to \infty} \left[ \tr\!\ V(\pi_N(r)) -\tr\!\ V(\pi_N(r^V)) - 2\sum_{k=1}^{N-1} \log (a_k/a^V_k) \right] = 
 \mathcal{K}(\mu_V\!\ |\!\ \mu) + \sum_{\lambda \in E} {\mathcal F}_{V}(\lambda) 
\end{align*} 
for $\mu\in \mathcal{S}_1(I)$ and if $\mu\notin \mathcal{S}_1(I)$, the left hand side equals $+\infty$. 
\end{thm}

This sum rule has to be compared with the sum rule proved by Nazarov et al. \cite{Nazarov}. 
Let us assume, without loss of generality, that we are in the one-cut regime with $I = [-2,2]$. 
Nazarov et al. 
 considered equilibrium measure of the form $A(x) \SC(dx)$ where $A$ is a nonegative polynomial. They define the function $F$ by
\begin{align}
F(x) = \begin{cases} \displaystyle\int_{2}^x A(t) \sqrt{(t-2)(t+2)} \diff t & x \geq 2\\
\displaystyle \int_x^{-2} A(t) \sqrt{(2-t) (t+2)} \diff t & x \leq -2\,,
\end{cases}
\end{align} 
and state the sum rule with $F$ instead of $\mathcal F_V$.

Let us recall how $A$ and $V$ are related. First,  it is known  (\cite{pastur_shcherbina} Th. 11.2.4) that % if $R(x) = (4-x^2)$, then 
if $V \in \mathcal C^2$, then 
\begin{align} \label{eq:polynomialmuV}
\mu_V(dx) = \frac{1}{2\pi}A(x)\sqrt{4-x^2} \diff x = A(x) \diff \SC(x)
\end{align}
with 
\begin{align} \label{VtoA}
A(x) = \frac{1}{\pi}\int_{-2}^2 \frac{V'(x) - V'(t)}{x-t} \frac{\diff t}{\sqrt{4-t^2}}\,.
\end{align}
In particular, if $V \in \mathcal V_{2d}$, then $A$ is a polynomial of degree $(2d-2)$.

Conversely, if $A \in \mathcal C^1$, we have
\begin{align}
\label{AtoV}
V'(x) = xA(x) -2 \int \frac{A(x) - A(t)}{x-t} \diff \SC(t)\,.
\end{align}
Indeed, let us consider 
 the master equation 
connecting $V$ and $\mu_V$ :
\begin{align*}\frac{V(x)}{2} - \int \log |x-y| d\mu_V (y) = \mathcal E(\mu_V) -\frac12 \int V d\mu_V \  \hbox{for}\ \  x \in [-2,2]\,.
\end{align*}
By differentiation, we get
\begin{align}
\label{master}
\frac{1}{2} V'(x) &= \mathrm{P.V.} \int \frac{\diff \mu_V(t)}{x-t}  \\ \notag
&=  \mathrm{P.V.} \int \frac{A(t)}{x-t}\diff \SC(t)\\ 
 \label{integrand} &=  \mathrm{P.V.}\left[\int \frac{A(t) - A(x)}{x-t}\diff \SC(t) + A(x) \int \frac{\diff \SC(t)}{x-t} \right]\\  \label{cont}
&= \int \frac{A(t) - A(x)}{x-t}\diff \SC(t) + A(x)\!\ \mathrm{P.V.}\int \frac{\diff \SC(t)}{x-t}
\\ \label{master2}&= \int \frac{A(t) - A(x)}{x-t}\diff \SC(t) + \frac{x}{2} A(x)\,, 
\end{align}
where (\ref{cont}) holds by continuity of the integrand in the first integral of (\ref{integrand}), and (\ref{master2}) holds by application of  (\ref{master}) to the potential $x^2/2$. This gives exactly (\ref{AtoV}) for $x \in [-2, 2]$, and 
 %We conclude that for $x \in [-2, 2]$
%\begin{align}
%\label{id}
%V'(x) = xA(x) -2 \int \frac{A(x) - A(t)}{x-t} \SC(dt)
%\end{align}
%\section{Randomization}
this can be extended for every $x$ real, 
 as an equality between polynomials.

As a consequence, from (\ref{poteff}) and (\ref{defF}) we have for $x\notin [-2, 2]$
\begin{align}
\mathcal F'_V (x) &= V'(x) - 2 \int \frac{A(t)}{x-t} \diff \SC(t) \notag \\
&= xA(x) - 2 A(x) \int \frac{\diff \SC(t)}{x-t}\notag \\
&= xA(x) - A(x) \left(x - \sqrt{x^2-4}\right)\notag \\
&= A(x) \sqrt{x^2-4}
\end{align}
hence (A2) is satisfied if $A \geq 0$.

Actually, it is known that, if $A > 0$ on $[-2,2]$ and if  (A2) is  satisfied, then 
$\mathcal F_V = F$ (a consequence of (1.13) in \cite{Albeverio}). In this case, (\ref{Lieb}) holds true.

\section{Large deviation results}
\label{sec:LDPs}

In order to be self-contained, we recall the basic definition and tools in the appendix, as well as a technical result used in the proof. We refer to \cite{demboz98} for more details.  
The classical LDP for the empirical eigenvalue measure (defined in \eqref{empiricallaw}) is widely known, see \cite{arous1997large} or \cite{agz}, Theorem 2.6.1. It holds in the space $\mathcal{M}_1$, equipped with the weak topology. 

\begin{thm} \label{thm:classicalLDP}
Assume that $V$ satisfies the assumption (A1). Then the sequence of empirical spectral measures $\mu_\u^{(n)}$ satisfies the LDP with speed $n^2$ and good rate function
\begin{align*}
\mathcal{I}_u (\mu) = \mathcal{E}(\mu) - \inf_\nu\mathcal{E}(\nu)
\end{align*}
where $\mathcal{E}$ is as in \eqref{ratemuu}. 
%\blue{$\mathcal E(\mu)$ achives its minimumu value uniquely in $\mu_V$.}
\end{thm}
This
 theorem shows that  $\mu_\u\sn$ %the empirical spectral measure
 converges somehow quickly towards the unique minimizer $\mu_V$ of $\mathcal I_u$ %the equilibrium measure
 since the speed in the LDP is $n^2$. On the other hand, %Now, if one is interested with
 the convergence of the extremal eigenvalue to the extremal point of the support of $\mu_V$ %the convergence
 is slower. Indeed,
the extremal eigenvalue satisfies an LDP only at speed $n$ as stated in the following theorem (see \cite{BorGui2013onecut} Prop. 2.1 based on \cite{arous2001aging}).

\begin{thm} \label{thm:classicalextremalLDP}
If $V$ is continuous and satisfies (A1) and (A3), then the random variable $\lambda\sn_{\max} = \max\{\lambda_1, \dots, \lambda_n\}$ satisfies the LDP at speed $n$ and good rate function $\mathcal F_V$.
\end{thm}

In the multi-cut case, we can also show an LDP for eigenvalues between two intervals of $I$ (Theorem \ref{thm:extremalLDP}). A related statement about %The statement for a collection of
 outliers is given in \cite[Lemma 3.1]{BorGuimulticut}.

\subsection{Spectral LDP}

The following theorem is our main large deviation result. The spectral measures $\mu_n$ are considered as random elements in $\mathcal M_1$, 
 equipped with the weak topology and the corresponding Borel $\sigma$-algebra.

\begin{thm} \label{MAIN}
Assume that the potential $V$ satisfies the assumptions (A1), (A2) and (A3). Then the sequence of spectral measures $\mu_n$ satisfies under $\Pnv$ the LDP with speed $\beta'n$ and good rate function
\begin{align*}
\mathcal{I}_{\operatorname{sp}}(\mu) = \mathcal{K}(\mu_V\!\ |\!\ \mu) + \sum_{\lambda \in E(\mu)} {\mathcal F}_{V}(\lambda)
\end{align*}
if $\mu \in \mathcal{S}_1(I)$ and $\mathcal{I}_V(\mu) = \infty$ otherwise. 
\end{thm}

\subsection{Coefficient LDP}

To obtain an expression for the rate function of the random recursion coefficients, we need to assume that $V$ is a polynomial. Recall that $\mathcal{M}_{1,K}$ is the set of probability measures with support in $[-K,K]$ and $\psi$ maps a measure $\mu$ to its Jacobi coefficients. The large deviation principle for the recursion coefficients will be under conditioning on the set $\mathcal R_K=\psi(\mathcal{M}_{1,K})$ and we define $\mathbb{P}^V_{n,K} = \mathbb{P}^V_{n}(\cdot |\mathcal R_K)$. Note that if $\mu$ is a spectral measure with support in $[-K,K]$, the Jacobi coefficients satisfy $|b_k|,|a_k|\leq K$, see e.g. Proposition 1.3.8 in \cite{simon2}.

Recall that $r^V$ is the sequence of Jacobi parameters of the equilibrium measure $\mu_V$ and we will always choose $K$ so large that the support of $\mu_V$ is contained in the interior of $\mathcal{M}_{1,K}$.

\begin{thm} \label{MAIN2}
Suppose $V$ is a polynomial of even degree $V$ with positive leading coefficient. Then the sequence $(\mu_n)_n$ satisfies under $\mathbb{P}^V_{n,K}$ the LDP in $\mathcal{M}_{1,K}$ with speed $n\beta'$ and good rate function
\begin{align*}
\mathcal{I}_{\operatorname{co},K}(\mu) = \lim_{N\to \infty} \left[ \tr V(\pi_N(r)) -\tr V(\pi_N(r^V)) - 2\sum_{k=1}^{N-1} \log (a_k/a^V_k) + \xi_{N,K}(\pi_N(r))\right] .
\end{align*}
where the term $\xi_{N,K}(\pi_N(r))$ satisfies
\begin{align}\label{crucialbound}
|\xi_{N,K}(\pi_N(r))| \leq C(K,V) \limsup_{N\to \infty}  M_+ (r_N)
%\big( |a_{J-d}-a^V_{J-d}| + |b_{J-d}-b^V_{J-d}| + \dots +  |b_{J+d}-b^V_{J+d}|\big)
\,,
\end{align} 
for some constant $C(K,V) > 0$ and
\begin{align}
\label{defM+}
M_+ (r_N) =  |a_{N-d}-a^V_{N-d}| + |b_{N-d}-b^V_{N-d}| + \dots +  |b_{N+d}-b^V_{N+d}|\,.
\end{align}

\end{thm}

The proof of this LDP is done in Section \ref{sec:coeffproof}. Actually 
 this proof is not independent from 
the proof of Theorem \ref{MAIN} since 
  it uses several times that we know 
 an LDP holds. However, the method of proof is different and uses directly the density \eqref{krishnadensity}.  Unfortunately, this does not give an explicit expression for the term $\xi_{N,K}$ in the rate function. The bound \eqref{crucialbound} %in Theorem \ref{MAIN2}
 implies though that this term is uniformly bounded on $\mathcal{R}_K$. In particular, it does not influence the finiteness of the rate function, which is crucial in view of the {\it gem} in Section \ref{sec:sumrules}.

We remark that although the density \eqref{krishnadensity} always gives a (finite) vector of Jacobi parameters, the above rate function is only finite if $a_k>0$ for all $k$, that is, $r$ is a Jacobi sequence of a measure with infinite support.

For general polynomial potentials $V$, the expression for the rate cannot be extended to the full space $\mathcal{R}$, as the constant $C(K,V)$ would blow up as $K\to \infty$. We do have a good control in the one-cut regime, where we can show that $\xi_K(r)$ does not change the value of the $K$-independent part of the rate. The consequence is the following LDP for the unrestricted measure.

\begin{thm} \label{MAIN3}
Suppose that the assumptions of Theorem \ref{MAIN2} hold. Assume further that the support of $\mu_V$ is a single interval. Then, the sequence $(\mu_n)_n$ satisfies under $\mathbb{P}^V_{n}$ the LDP with speed $n\beta'$ and good rate function
\begin{align*}
\mathcal{I}_{\operatorname{co}}(\mu) = \lim_{N\to \infty} \left[ \tr V(\pi_N(r)) -\tr V(\pi_N(r^V)) - 2\sum_{k=1}^{N-1} \log (a_k/a^V_k) \right] .
\end{align*}
\end{thm}

\medskip

\begin{rem} \label{rem:errorterm}
It is an interesting question to ask for the role of the remainder term $\xi_{N,K}$ in Theorem \ref{MAIN2}. Considering that the one-cut-LDP in Theorem \ref{MAIN3} does not involve such a term, either it is an artifact or an immanent feature of the multi-cut-case. We argue that the remainder is in fact not an artifact, but necessary in order to distinguish between different measures obtained by permuting the Jacobi coefficients. 

As an example, we consider the quartic potential $V(x) = \frac{x^4}{4}-\frac{\mathrm{v}x^2}{2}$. When $\mathrm{v}>2$, the equilibrium measure $\mu_V$ is supported by two disjoint intervals $[-\alpha^+, -\alpha^-] \cup [\alpha^-,\alpha^+]$ with $\alpha^\pm = \sqrt{\mathrm{v}\pm 2}$, see \cite{pastur_shcherbina}, Example 11.2.11 (2) and Problem 11.4.13 or \cite{Blower}, Section 4.6.
The Jacobi coefficients of $\mu_V$ are given by $b_k=0$ for all $k\geq 1$ and the $a_k$ are perturbations of two-periodic coefficients, as given in equation (14.2.16) in \cite{pastur_shcherbina}. Indeed, the measure $\mu_V$ is symmetric and the monic orthogonal polynomials of degree $2n$ may be written as $P_n(x^2)$, where the $P_n$ are orthogonal with respect to a measure $\mu_{\mathrm{even}}$ supported by $[\v-2,\v+2]$ and satisfy the recursion
\begin{align}
x P_n(x) = P_{n+1}(x) + (a_{2n}^2 + a_{2n+1}^2)  P_n(x) + a_{2n}^2a_{2n-1}^2 P_{2n-1}(x) .
\end{align}
Additionally, the monic orthogonal poylnomials of $\mu_V$ of degree $2n+1$ may be written as $xQ_n(x^2)$, with the $Q_n$ orthogonal to a measure $\mu_{\mathrm{odd}}$ supported by $[\v-2,\v+2]$, and they satisfy
\begin{align}
x Q_n(x) = Q_{n+1}(x) + (a_{2n+2}^2 + a_{2n+1}^2)  Q_n(x) + a_{2n+1}^2a_{2n}^2 Q_{2n-1}(x) .
\end{align}
The combination of the two recursions implies that 
\begin{align}
\label{Rakh}
\lim_{n\to \infty} a_{n}^2a_{n-1}^2 =1,  \qquad \lim_{n\to\infty} (a_{n}^2 + a_{n-1}^2) = \v .
\end{align}
This shows that (at least along subsequential limits, which we may ignore for the following argument) 
$a_{2n-1}\to a$, $a_{2n}\to \bar a$, where $a,\bar a$ are the two solutions to $\ell^2 - \v \ell +1 = 0$, i.e., 
\begin{align}
\ell_1 = \frac{\v- \sqrt{\v^2 -4}}{2}, \qquad \ell_1 = \frac{\v+ \sqrt{\v^2 -4}}{2} . 
\end{align}

Switching the $a_k$ of even and odd index, we obtain a new measure $\bar \mu$. If there would be no remainder term $\xi_{J,K}$, the rate function at $\bar \mu$ would be the limit as $N\to \infty $ of 
\begin{align} \label{eq:pseudorate}
\tilde{\mathcal{I}}^{(N)}_{\operatorname{co},K}(\bar \mu) =  \tr V(\pi_N(r)) -\tr V(\pi_N(r^V)) - 2\sum_{k=1}^{N-1} \log (a_k/a^V_k) .
\end{align} 
However, the quasi-periodic structure of $\mu$ causes an alternating behavior of $\tilde{\mathcal{I}}^{(N)}_{\operatorname{co},K}(\bar \mu)$. 
More precisely, a straightforward but lengthy calculation yields that 
\begin{align} \label{eq:pseudorate2}
\lim_{N\to \infty} \left( \tilde{\mathcal{I}}^{(2N)}_{\operatorname{co},K}(\bar \mu)-\tilde{\mathcal{I}}^{(2N-1)}_{\operatorname{co},K}(\bar \mu)\right)
= \frac{1}{2}({\bar a}^4-a^4) - \mathrm{v}({\bar a}^2-a^2)-2\log(\bar a/a) , 
\end{align}
so that \eqref{eq:pseudorate} does not converge as $N\to \infty$. 
\end{rem} 

\medskip

\subsection{From LDPs to sum rules}

In this section, we prove Theorem \ref{abstractgem} and Theorem \ref{newsumrule}. The main argument in both cases is that we have two different expressions for the large deviation rate function, one using the spectral encoding and one using the encoding
 by Jacobi coefficients. Since both expressions must agree, they yield the ``spectral side'' and the ``coefficient side'', respectively, of a sum rule. 

\medskip

\textbf{Proof of Theorem \ref{newsumrule}:} Suppose $V$ is a nonzero polynomial of even degree, such that the equilibrium measure $\mu_V$ is supported by a single interval $I$. Theorem \ref{MAIN} yields the LDP for $(\mu_n)_n$ with speed $n$ and rate $\mathcal{I}_{\operatorname{sp}}$. On the other hand, by Theorem \ref{MAIN3}, $(\mu_n)_n$ satisfies the LDP with speed $n$ and rate function $\mathcal{I}_{\operatorname{co}}$. Since a large deviation rate function is unique, we have the equality
\begin{align}
\mathcal{I}_{\operatorname{sp}}(\mu) = \mathcal{I}_{\operatorname{co}}(\mu)
\end{align}
for all $\mu \in \mathcal{M}_1$. For $\mu \in \mathcal{S}_1(I)$, this is precisely the equality claimed in Theorem \ref{newsumrule}. For $\mu \notin \mathcal{S}_1(I)$, we know that the left hand side satisfies $\mathcal{I}_{\operatorname{sp}}(\mu)=+\infty$, so the right hand side must equal $+\infty$ as well. \hfill $\Box$  

\medskip

\textbf{Proof of Theorem \ref{abstractgem}:} Let $V$ be a nonzero polynomial of even degree. We want to combine the LDP results of Theorem \ref{MAIN} and Theorem \ref{MAIN2}. The former are obtained under $\Pnv$, whereas the latter are under $\mathbb{P}^V_{n,K}$, for $K$ large enough (depending on $V$). Theorem \ref{thm:restrictedLDP} shows that $(\mu_n)_n$ satisfies also the LDP under $\mathbb{P}^V_{n,K}$ in the restricted space $\mathcal{M}_{1,K}$, with rate $\mathcal{I}_{\operatorname{sp},K}$ the restriction of $\mathcal{I}_{\operatorname{sp}}$ to $\mathcal{M}_{1,K}$. By uniqueness of rate functions, we obtain the restricted sum rule
\begin{align}
\mathcal{I}_{\operatorname{sp},K}(\mu) = \mathcal{I}_{\operatorname{co},K}(\mu)
\end{align}
for any $\mu \in \mathcal{M}_{1,K}$. For $\mu$ a probability measure with compact, infinite support, we may choose $K$ so large that $\mu \in \mathcal{M}_{1,K}$. Then the above equality holds, where both sides are simultaneously finite or infinite. Condition \eqref{gengem} is equivalent to finiteness of $\mathcal{I}_{\operatorname{co},K}(\mu)$ since from (\ref{crucialbound}) and (\ref{defM+})
\[ |\xi_{N,K}(\pi_N (r))| \leq 4(2d+1) K C(K,V)\,.\]By the restricted sum rule, this is equivalent to finiteness of  $\mathcal{I}_{\operatorname{sp},K}(\mu)$. We have $\mathcal{I}_{\operatorname{sp},K}(\mu)<\infty$ if and only if $\mu \in \mathcal{S}_1(I)$, and both $\sum_{\lambda\in E}\mathcal{F}_V(\lambda)$ and $\mathcal{K}(\mu_V|\mu)$ are finite. The first two conditions are just \emph{1.} and \emph{2.} in Theorem \ref{abstractgem}, and the third one is equivalent to \emph{3.}, see \eqref{eq:KL2}.
\hfill $\Box$

\section{Proof of the spectral LDP}
\label{sec:spectralproof}

\subsection{Structure of the proof}

This section is devoted to the proof of Theorem \ref{MAIN}. 
In the large deviation behavior of the weighted spectral measure $\mu_n$, all  eigenvalues outside of $I$ (the outliers) will contribute %be visible
 and in fact, the rate function in Theorem \ref{MAIN} can be finite even for countably many outliers. 
The main difficulty of the proof comes from 
%Checking
 the {\it a priori} unbounded number of eigenvalues close to $\partial I$, 
%the boundary of the intervals $I$,
and the dependence with the bulk of eigenvalues on $I$. 
As in our proof in the one-cut regime, the main idea is to apply the projective limit method to reduce the spectral measure to a measure with only a fixed number of eigenvalues outside the limit support $I$.
However, controlling the eigenvalues between two intervals in $I$ requires special care. We do this by dividing the outliers into groups according to which of the sub-intervals constituting $I$ they are the closest. This allows to apply the general strategy of the one-cut case, albeit in a much more technical way. Additionally, our new encoding for spectral measures introduced below also takes care of topological problems which occurred in \cite{magicrules} when transferring the LDPs from one space to another.

The main steps of the proof are as follows. To begin with, we decouple the weights of the random measure $\mu_n$ and introduce a non-normalized random measure $\tilde\mu_n$ with weights from a family of independent random variables. In the next section, we will introduce a family $\zeta(\tilde\mu_n)$ of points not in $\operatorname{Int}(I)$ encoding the outlying support points and a family $\gamma(\mu_n)$ of corresponding weights. If $\tilde\mu_{I,n}$ denotes the restriction of $\tilde\mu_n$ to $I$, then we may identify $\tilde\mu_n$ with the collection 
\begin{align*}
\big( \tilde\mu_{I,n}, \zeta(\tilde\mu_n), \gamma(\tilde\mu_n) \big) .
\end{align*} 
The LDP for $\mu_n$ is then proved using this representation, with the following intermediate steps. 
\begin{itemize}
\item[(1)] We prove an LDP for a finite collection of entries of $\zeta(\tilde\mu_n)$. This is Theorem \ref{thm:extremalLDP}. 
\item[(2)] Using (conditional) independence of the outliers $\zeta(\tilde\mu_n)$ and the weights $\gamma(\tilde\mu_n)$, we can prove in Theorem \ref{thm:jointLDP} a joint LDP for a finite collection of entries of $(\zeta(\tilde\mu_n),\gamma(\tilde\mu_n))$.   
\item[(3)] In Theorem \ref{thm:LDPprojectivelimit} we use the projective method (the Dawson-G\"artner Theorem) to prove the LDP for the whole family $(\zeta(\tilde\mu_n),\gamma(\tilde\mu_n))$.
\item[(4)] Once we have this LDP, we can use the technical result of Theorem \ref{newgeneral} to combine the outliers with $\tilde\mu_{I,n}$, and we show in Theorem \ref{thm:jointjointLDP} the joint LDP for $(\tilde\mu_{I,n},\zeta(\tilde\mu_n),\gamma(\tilde\mu_n))_n$.
\item[(5)] Finally, in Section \ref{sec:normalizing}, we use the contraction principle to transfer this LDP to the non-normalized spectral measure $\tilde\mu_n$ and recover after normalizing the spectral measure $\mu_n$. 
\end{itemize}

\subsection{New encoding of measures}
\label{sec:newencoding}

Let $\mu \in \mathcal{S}(I)$ be a nonnegative measure with the restrictions on the support as in \eqref{support}. Then, $\mu$ can be written as 
\begin{align}\label{muinS}
\mu = \mu_{I} +  \sum_{\lambda \in E(\mu)} \gamma_\lambda \delta_{\lambda} ,
\end{align}
where $\mu_I$ is the restriction to $I$. 
We now introduce a particular enumeration of the elements of $E(\mu)$, according to the %boundary
 point of $\partial I$ to which they are closest, and then according to their distance to that point. 
For this, recall that $I$ is a disjoint union of compact intervals $I_m$, and suppose $I_m=[l_m,r_m]$, so that $r_m<l_{m+1}$ for $m=1,\dots,M-1$. Let $\theta_m = \tfrac{1}{2}(r_m+l_{m+1})$ denote the midpoint between $I_{m}$ and $I_{m+1}$. Then there is a unique array $ \zeta = (\zeta_{i,j})_{i,j}$ $i=1,\dots 2M$ and $j\geq 1$ encoding the elements of $E(\mu)$, which is defined as follows :
 \begin{itemize}
\item
 $\zeta_{1,j}$ for $j\geq 1 $ are the elements of $E(\mu)$ to the left of $l_1$, in increasing order, 
\item
$\zeta_{2,j}$ is are the elements of $E(\mu)$ in $(r_1,\theta_1]$ in decreasing order, 
\item $\zeta_{3,j}$ are are in increasing order the elements in $E(\mu)\cap (\theta_1,l_2)$, 
\item and so on.
\end{itemize}  If there are only a finite number of such elements, the sequence $\zeta_{i,j}$ is extended by the boundary element $l_m$ for $i=2m-1$ and by $r_m$ for $i=2m$.  
More precisely, given $\mu\in \mathcal{S}(I)$, let $ \zeta = \zeta(\mu) = (\zeta_{i,j})_{i,j}$ be the unique array, such that  
\begin{align} \label{newencoding}
E(\mu)  =  \bigcup_{i=1}^{2M}  \bigcup_{j= 1}^\infty  \{\zeta_{i,j}\} \setminus  \partial I , 
\end{align}
and additionally, for $j \geq 1$,
\begin{align} \label{newencoding1}
\begin{split} \zeta_{1,j} \in (-\infty,l_1], & \qquad  \zeta_{2M,j} \in [r_M,\infty), \qquad  \\  
 \zeta_{2m,j}\in [r_m,\theta_m], & \qquad  \zeta_{2m+1,j}\in(\theta_m,l_{m+1}] , 
\end{split}
\end{align}
for $m=1,\dots , M-1$, and for all $i\leq 2M, j \geq 1$,
\begin{align} \label{newencoding2}
d(\zeta_{i,j},I)\geq d(\zeta_{i,j+1},I) \quad \text{  and }\quad  d(\zeta_{i,j},I)> d(\zeta_{i,j+1},I)\text{  unless } \zeta_{i,j}\in \partial I
\end{align}
(recall that $d(\cdot ,I)$ %measures
 is the distance to the set $I$).  
Condition \eqref{newencoding1} ensures that the elements are grouped according to the closest point in $\partial I$, and condition \eqref{newencoding2} ensures that the elements are strictly ordered, unless there are only finitely many. 
The union of all entries in $\zeta$ as in \eqref{newencoding} yields the elements of $E$ again, and in addition possibly the boundary points if there are only finitely many nonzero entries. We denote the closure (in the product topology on $\mathbb{R}^{2M\times \mathbb{N}}$) of the set of all arrays $(\zeta_{i,j})_{i,j}$ satisfying \eqref{newencoding1} and \eqref{newencoding2} by $\mathcal{Z}$. 
%Then $\mathcal{Z}$ is a closed subset of a set compact by Tychonoff's Theorem, and therefore compact as well. 

In order to encode the weights of a measure as in \eqref{muinS} as well, let $\gamma(\mu) =(\gamma_{i,j})_{i,j}$, $i=1,\dots ,2M, j\geq 1$ be the unique non-negative array such that
\begin{align} \label{newencoding3}
\gamma_{i,j} = 0 \quad \text{ if and only if }\quad  \zeta_{i,j} \in \partial I
\end{align}
and such that
\begin{align} \label{newencoding4}
\mu = \mu_{I} + \sum_{i=1}^{2M}\sum_{j=1}^\infty \gamma_{i,j} \delta_{\zeta_{i,j}}. 
\end{align}
The set of weights is denoted by 
\begin{align} \label{dafweightset}
\mathcal{G} =  [0,\infty)^{2M\times \mathbb{N}}
%\mathcal{G} = \left\{ \gamma \in [0,\infty)^{2M\times \mathbb{N}} \, \left| \, \sum_{i,j} \gamma_{i,j} %<\infty \right. \right\}  
\end{align} 
and we endow $\mathcal{G}$ with the product topology. 
These definitions set up a one-to-one correspondence between a finite measure $\mu \in \mathcal{S}(I)$ and 
\begin{align} \label{newencodingfinal}
\big( \mu_{I}, \zeta(\mu), \gamma(\mu) \big)  \in \mathcal{M}(I)\times \mathcal{Z}\times \mathcal{G} , 
\end{align}
where $(\zeta, \gamma)$ satisfy \eqref{newencoding3}.

The representation \eqref{newencodingfinal} will be applied not directly to the spectral measure $\mu_n$, but to a variant with uncoupled, independent weights. Recall that under $\Pnv$, the vector $(w_1,\dots ,w_n)$ is Dirichlet distributed and has the same distribution as 
\begin{align}\label{decoupling}
\left(\frac{\omega_1}{\omega_1+ \dots + \omega_n}, \dots , \frac{\omega_n}{\omega_1+ \dots + \omega_n}\right) ,
\end{align}
where $\omega_1, \dots, \omega_n$ are independent variables 
with distribution Gamma$(\beta', (\beta'n)^{-1})$ and mean $n^{-1}$. Without loss of generality, assume that the variables $\omega_k$ are defined on the same probability space as the $\lambda$'s and independent of them. 
We then consider the non-normalized measure 
\begin{align}\label{unnorm}
\tilde\mu_n = \sum_{k=1}^n \omega_k \delta_{\lambda_k} \in \Sr (I)
\end{align}
and we can come back to the original measure by normalization.
%ing $\tilde\mu_n(\mathbb{R})^{-1}\tilde\mu_n$. 
Therefore, we start by looking at 
\begin{align*}
\big( \tilde\mu_{n,I}, \zeta(\tilde\mu_n), \gamma(\tilde\mu_n) \big)
\end{align*}
and to simplify notation, we will write $\zeta^{(n)}$ for $\zeta(\tilde\mu_n)$ and $\gamma^{(n)}$ for $\gamma(\tilde\mu_n)$. 

\subsection{LDP for a finite collection of extremal eigenvalues}
\label{sec:spectralproof1}

In this section, we prove an LDP  %large deviation principle
 for a finite collection of elements of the arrays $(\zeta^{(n)})_n= (\zeta(\tilde\mu_n))_n$. Fix a $N \geq 1$ and let 
\begin{align} \label{defproject}
\pi_N:\mathbb{R}^{2M\times \mathbb{N}} \to \mathbb{R}^{2M\times N}
\end{align}
denote the canonical projection onto the first $N$ columns. We denote by $\zeta_N^{(n)}= \pi_N(\zeta^{(n)})$ the image of the outlying support points and let $\mathcal{Z}_N = \pi_N(\mathcal{Z})$. 
The following LDP for the finite collection of extremal eigenvalues is a crucial starting point for the LDP of $\mu_n$.

\begin{thm} \label{thm:extremalLDP}
Under $\mathbb{P}_n^V$, the collection of extreme eigenvalues $(\zeta_N^{(n)})_n$ satisfies the LDP in $\mathcal{Z}_N$ with speed $\beta' n$ and good rate function
\begin{align*}
\mathcal{I}_N^{\mathrm{ext}}(z) = 
\sum_{i=1}^{2M} \sum_{j=1}^N \mathcal{F}_V( z_{i,j}) .  
\end{align*}
\end{thm}

The proof follows the main steps of Theorem 4.1 in \cite{magicrules}. Therein,  the $N$ largest and smallest eigenvalues were considered. The multi-cut situation, besides being notationally heavier, requires some additional care. This is not only due to outliers between two intervals in $I$, but also to the new encoding of outliers. This encoding was not useful in the one-cut case. For this reason, we give the main arguments of the proof in Section \ref{sect:extremalLDP1}, and refer to \cite{magicrules} for the detailed calculations.

The next main step is then to combine the finite collection of extremal eigenvalues with their weights. 

\subsection{LDP for a finite collection of eigenvalues and weights}
\label{sec:spectralproof2}

Similarly to the definition in Section \ref{sec:spectralproof1}, we denote by $\gamma_N^{(n)}= \pi_N(\gamma^{(n)})$ the projection of the array of weights and let $\mathcal{G}_N = \pi_N(\mathcal{G})$. The following joint LDP for $\zeta_N^{(n)}$ and $\gamma_N^{(n)}$ is the main result in this section. Since $\zeta_{i,j}^{(n)}\in \partial I$ implies in our encoding that $\gamma_{i,j}^{(n)}=0$, the two arrays are not independent. However, conditioned on $\{\zeta_{i,j}^{(n)}\notin \partial I\}$, the eigenvalue $\zeta_{i,j}^{(n)}$ and its weight $\gamma_{i,j}^{(n)}$ are actually independent. Using this fact, and the explicit (conditional) distribution of $\gamma_{i,j}^{(n)}$, the proof becomes fairly straightforward. 

Let us remark that we prove the joint LDP in the ``full'' space $\mathcal{Z}_N \times \mathcal{G}_N$ without the above condition on some weights being zero, as formalized in \eqref{newencoding3}. While the distribution $\Pnv$ is concentrated on the subset satisfying \eqref{newencoding3}, this would lead to a rate function without compact level sets. In view of later parts of the proof, we consider the larger space with the lower semi-continuous continuation of the rate function.

\begin{thm} \label{thm:jointLDP}
For any $N\geq 1$, the sequence $(\zeta_N^{(n)},\gamma_N^{(n)})_n$ satisfies under  $\Pnv$ the LDP in $\mathcal{Z}_N \times \mathcal{G}_N$ with speed $\beta'n$ and good rate function
\begin{align*}
\mathcal{I}_N^{(\mathrm{ext},\mathrm{w})}(z,g) = \mathcal{I}_N^{\mathrm{ext}}(z) + ||g||_{N,1}, 
\end{align*}
with $||\cdot ||_{N,1}$ the $\ell_1$-norm on $\mathcal{G}_N$. 
\end{thm}

\textbf{Proof:} Let $\tilde \gamma^{(n)}_{i,j}$, $1\leq i \leq 2M, j\geq 1$ be independent and  Gamma$(\beta', (\beta'n)^{-1})$ distributed random variables, defined on the same probability space as $\zeta_N^{(n)},\gamma_N^{(n)}$, and independent of $\zeta_N(\tilde\mu_n)$. Then, by \eqref{unnorm}, we have the equality in distribution 
\begin{align} \label{equaldistextweights}
\big(\zeta^{(n)}_{i,j},\gamma_{i,j}^{(n)}\big)_{i,j} \stackrel{d}{=} \big( \zeta^{(n)}_{i,j},\tilde\gamma^{(n)}_{i,j}\mathbbm{1}_{\{ \zeta_{i,j}^{(n)}\notin \partial I \}} \big)_{i,j}. 
\end{align}
Let $\tilde{\gamma}^{(n)}_N = \pi_N( (\tilde{\gamma}^{(n)}_{i,j})_{i,j})$. It follows by straightforward calculations, that for each $i,j$, the sequence $(\tilde{\gamma}^{(n)}_{i,j})_n$ satisfies the LDP in $[0,\infty)$ with speed $\beta'n$ and good rate function $I_0$, with $I_0(x) = x$. Since the $\tilde\gamma_{i,j}$ are independent, this implies the LDP for $\tilde{\gamma}^{(n)}_N$ in $\mathcal{G}_N$ with speed $\beta'n$ and good rate function 
\begin{align}
\mathcal{I}_N^{\mathrm{w}}(\tilde{g}) = \sum_{i=1}^{2M} \sum_{j=1}^N I_0(\tilde{g}_{i,j}) = ||\tilde{g}||_{N,1} . 
\end{align}
For the joint LDP, let us consider the finite family $(\zeta_{1,j}^{(n)},\gamma_{1,j}^{(n)})_{1\leq j\leq N}$ corresponding to eigenvalues to the left of the leftmost interval. Then, for sets $A=A_1\times \dots \times A_N, B=B_1\times \dots\times B_N\subset \mathbb{R}^N$, \eqref{equaldistextweights} implies
\begin{align} \label{jointLDP0}
& \Pnv\big((\zeta_{1,j}^{(n)})_j \in A,(\gamma_{1,j}^{(n)})_j \in B )
  = \Pnv\big((\zeta_{1,j}^{(n)})_j \in A )\Pnv( (\tilde\gamma_{1,j}^{(n)})_j \in B) ,
\end{align}  
whenever $\partial I \cap A_N= \emptyset$, that is, we require the rightmost outlier (and then all of them) to be outside of $I$. The LDPs for $((\zeta_{1,j}^{(n)})_j)_n$ and for $((\tilde\gamma_{1,j}^{(n)})_j)_n$ (which can be obtained from the LDPs for $(\zeta_{N}^{(n)})_n$ and for $(\tilde\gamma_{N}^{(n)})_n$ by the contraction principle), imply then for any $A,B$ as above and closed
\begin{align} \label{jointLDPub}
\limsup_{n\to \infty} \frac{1}{\beta'n}\log \Pnv\big((\zeta_{1,j}^{(n)})_j \in A,(\gamma_{1,j}^{(n)})_j \in B ) 
& \leq - \inf_{z \in A\cap \mathcal{Z}_N} \mathcal{I}_N^{\mathrm{ext}}(z) - \inf_{g \in B\cap \mathcal{G}_N} \mathcal{I}_N^{\mathrm{w}}(g) \notag \\
& = - \inf_{(z,g) \in (A \times B)\cap (\mathcal{Z}_N\times \mathcal{G}_N)} \left( \mathcal{I}_N^{\mathrm{ext}}(z) + ||g||_{N,1}\right) .   
\end{align}  
For $A,B$ as above and open, we get the lower bound
\begin{align} \label{jointLDPlb}
\liminf_{n\to \infty} \frac{1}{\beta'n}\log \Pnv\big((\zeta_{1,j}^{(n)})_j \in A,(\gamma_{1,j}^{(n)})_j \in B ) 
& \geq - \inf_{(z,g) \in (A \times B)\cap (\mathcal{Z}_N\times \mathcal{G}_N)} \left( \mathcal{I}_N^{\mathrm{ext}}(z) + ||g||_{N,1}\right) .   
\end{align}  
In fact, the lower bound can easily be extended to open sets $A=A_1\times \dots \times A_N, B=B_1\times \dots\times B_N$, with $A_j,B_j$ generic open subsets of $\mathbb{R}$. Set $A' = A\setminus I^N$, then 
\begin{align}
\Pnv\big((\zeta_{1,j}^{(n)})_j \in A,(\gamma_{1,j}^{(n)})_j \in B ) \geq 
\Pnv\big((\zeta_{1,j}^{(n)})_j \in A',(\gamma_{1,j}^{(n)})_j \in B ) , 
\end{align} 
and since $A'$ is still an open set, the generic lower bound follows from \eqref{jointLDPlb}. For the general upper bound, let $A,B$ be again of product form as above and closed. We define a modification $B' = B_1'\times \dots \times B_N'$ as follows. If $\partial I \cap A_j = \emptyset$, or $0\notin B_j$, set $B_j'=B_j$. If, on the other hand, $\partial I \cap A_j \neq \emptyset$, and $0\in B_j$, set $B_j'=[0,\infty)$. Then we have 
\begin{align} \label{jointLDPub2}
\Pnv\big(\zeta_{1,j} \in A_j,\gamma_{1,j} \in B_j ) \leq
\Pnv\big(\zeta_{1,j} \in A_j,\tilde \gamma_{1,j} \in B_j' ) =
\Pnv\big(\zeta_{1,j} \in A_j)\Pnv(\tilde \gamma_{1,j} \in B_j' ) .
\end{align}
This extends also to the whole vector, yielding 
\begin{align} \label{jointLDPub3}
\limsup_{n\to \infty} \frac{1}{\beta'n}\log \Pnv\big((\zeta_{1,j}^{(n)})_j \in A,(\gamma_{1,j}^{(n)})_j \in B ) 
 \leq - \inf_{(z,g) \in A \times B'\cap \mathcal{Z}_N\times \mathcal{G}_N} \left( \mathcal{I}_N^{\mathrm{ext}}(z) + ||g||_{N,1}\right) .   
\end{align} 
The general upper bound follows then from \eqref{jointLDPub3}, since the infimum over $g\in B'$ may be replaced by the infimum over $g\in B$. This implies the LDP for $(\zeta_{1,j}^{(n)},\gamma_{1,j}^{(n)})_{1\leq j\leq N}$. The arguments can be directly extended to outliers and weights in each of the intervals, which yields then the LDP for the family $(\zeta_N^{(n)},\gamma_N^{(n)})_n$. 
\hfill $\Box$

\subsection{LDP for the projective limit of extremal eigenvalues and weights}

By Theorem \ref{thm:jointLDP}, each projected sequence $(\zeta_N^{(n)},\gamma_N^{(n)})_n$ satisfies an LDP with a good rate function. We can then apply the Dawson-G\"artner Theorem (see the Appendix). It yields the LDP for the sequence of projective limits
\begin{align}
(\pi_N(\zeta^{(n)}),\pi_N(\gamma^{(n)}))_{N\geq 1}
\end{align}
in the projective limit of the spaces $\pi_N(\mathcal{Z})\times \pi_N(\mathcal{G})$. Since the topology on $\mathcal{Z} \times \mathcal{G}$ is the product topology, the canonical embedding from the projective limit into $\mathcal{Z} \times \mathcal{G}$ is continuous. An application of the contraction principle yield then the following result.

\begin{thm} \label{thm:LDPprojectivelimit}
The sequence $(\zeta^{(n)},\gamma^{(n)})_n$ satisfies under $\Pnv$ the LDP in $\mathcal{Z}\times \mathcal{G}$ with speed $\beta'n$ and good rate function
\begin{align*}
\mathcal{I}^{(\mathrm{ext},\mathrm{w})}(z,g)= \sup_{N\geq 1} \mathcal{I}_N^{(\mathrm{ext},\mathrm{w})}(z,g) = \sum_{i=1}^{2M} \sum_{j=1}^\infty \mathcal{F}_V( z_{i,j}) + |g_{i,j}| . 
\end{align*}
\end{thm}

\subsection{Joint LDP for the measure on $I$, the extremal eigenvalues and the weights}

The main result in this subsection is the following joint LDP, when we also add $\tilde\mu_{n,I}$, the restriction of $\tilde\mu_n$ to $I$.

\begin{thm} \label{thm:jointjointLDP}
The sequence $(\tilde\mu_{n,I},\zeta^{(n)},\gamma^{(n)})_n$ satisfies under $\Pnv$ the LDP in $\mathcal{M}(I)\times \mathcal{Z}\times \mathcal{G}$ with speed $\beta'n$ and good rate function
\begin{align*}
\widetilde{\mathcal{I}}(\tilde\mu,z,g) = \mathcal{K}(\mu_V | {\tilde\mu}) + \tilde\mu(I)-1+\mathcal{I}^{(\mathrm{ext},\mathrm{w})}(z,g) . 
\end{align*}
\end{thm}

\textbf{Proof:} 
The proof makes use of the LDP in Theorem \ref{thm:LDPprojectivelimit} for the extremal eigenvalues and their weights, and Theorem \ref{newgeneral} to combine this with the measure restricted to $I$. 

We check the conditions of Theorem \ref{newgeneral}, beginning with exponential tightness. The set 
\begin{align}
K_{H,T} = \big\{ (\mu,z,g)\in \mathcal{M}(I)\times \mathcal{Z}\times \mathcal{G}\, \big| \, 
||z||_\infty \leq H, \mu(I)+||g||_1 \leq T \big\}
\end{align}
is compact, and for $H$ so large that $I\subset [-H,H]$, 
\begin{align} \label{exptight}
\Pnv\left( (\tilde\mu_{n,I},\zeta^{(n)},\gamma^{(n)}) \notin K_{H,T} \right)
& \leq \Pnv\left( \zeta_{1,1}<-H\right) + \Pnv\left( \zeta_{2M,1}>H) \right) \notag \\
& \qquad + \Pnv\left( \omega_1+\dots + \omega_n> T\right) .
\end{align}
By Theorem \ref{thm:extremalLDP} (LDP for extremal eigenvalues), we have
\begin{align} \label{exptight2}
\limsup_{n\to \infty} \frac{1}{\beta' n} \log \Pnv\left( \zeta_{1,1}<-H\right)
\leq - \inf_{x\leq H} \mathcal{F}_V(x) , 
\end{align}
From the definition of $\mathcal{F}_V$ in \eqref{defF} we see that the upper bound goes to $-\infty$ as $H\to \infty$. For the last probability in \eqref{exptight}, we have by Cram\'er's Theorem for Gamma-distributed random variables,
\begin{align}\label{exptight3}
\limsup_{n\to \infty} \frac{1}{\beta' n} \log \Pnv\left( \omega_1+\dots + \omega_n> T\right) \leq - \big( T- \log T -1 \big) .  
\end{align}
Combining \eqref{exptight2} (and the analogous bound for the largest eigenvalue) and \eqref{exptight3}, we see that
\begin{align}\label{exptight4}
\lim_{H,T \to \infty} \limsup_{n\to \infty} \frac{1}{\beta' n} \log 
\Pnv\left( (\tilde\mu_{n,I},\zeta^{(n)},\gamma^{(n)}) \notin K_{H,T} \right) = -\infty ,
\end{align}
that is, the sequence  $(\tilde\mu_{n,I},\zeta^{(n)},\gamma^{(n)})_n$ is exponentially tight.

Now, let $D$ be the set of continuous $f:I\to (-\infty,1)$ and let $\varphi\in C_b(\mathcal{Z}\times \mathcal{G})$. We need to calculate the limit, on a logarithmic scale, of 
\begin{align} \label{jointmgf}
\mathcal{G}_n(f,\varphi) :&= \mathbb E_n^V \left[ \exp \left(n\beta' \int f\, d \tilde \mu_{n,I} + n\beta' \varphi( \zeta^{(n)},\gamma^{(n)} )\right)\right] \notag \\
&= \mathbb E_n^V \left[ \exp \left(n\beta'\sum_{k:\lambda_k\in I} \omega_k f(\lambda_k)  + n\beta' \varphi(\zeta^{(n)},\gamma^{(n)} )\right)\right] .
\end{align}
We will see that the main reasons which allow us to calculate the limit is the independence of the decoupled weights and then the faster LDP for the sequence of empirical spectral measures $\mu_n$. 
Indeed, recall that the weights $\omega_1,\dots ,\omega_n$ are independent and  Gamma$(\beta', (\beta'n)^{-1})$ distributed and, conditioned on the eigenvalues $\lambda_1,\dots ,\lambda_n$, the weights $(\omega_k)_{ \lambda_k \in I }$ are independent of $\zeta^{(n)}$. For each individual weight $\omega_k$ we have 
\begin{align*}
\frac{1}{\beta'} \mathbb{E}_n^V [e^{n\beta' t \omega_k}] = L(t). 
\end{align*}
Conditioning in \eqref{jointmgf} on $\lambda$ and integrating with respect to $(\omega_k)_{ \lambda_k \in I }$ yields therefore
\begin{align}
\mathcal{G}_n(f,\varphi) &= \mathbb E_n^V \left[ \mathbb E_n \left[ \left. \exp \left(n\beta'\sum_{k:\lambda_k\in I} \omega_k f(\lambda_k)  + n\beta' \varphi(\zeta{(n)},\gamma^{(n)} )\right) \right|\, \lambda \right]\right] \notag \\
&= \mathbb E_n^V \left[ \exp \left(n\beta'\int (L\circ f)\, d\mu_{n,I}^{(\u)}\right) \mathbb E_n \left[ \left. \exp \left( n\beta' \varphi( \zeta^{(n)},\gamma^{(n)} )\right)\right|\, \lambda \right]\right] \notag \\
&= \mathbb E_n^V \left[ \exp \left(n\beta'\int (L\circ f)\, d\mu_{n,I}^{(\u)} +  n\beta' \varphi( \zeta^{(n)},\gamma^{(n)} )\right) \right] ,
\end{align}
where $\muunI $ is the restriction of $\muun$ to $I$. We may now proceed as in \cite{magicrules}, Section 4.2. The empirical eigenvalue measure $\muun$ (and then also the restriction $\muunI$) satisfies the LDP at the faster scale $n^2$, which allows to replace it at our slower scale by its limit $\mu_V$. This yields
\begin{align}
\lim_{n\to \infty} \frac{1}{\beta'n} \log \mathcal{G}_n(f,\varphi) & = \lim_{n\to \infty} \frac{1}{\beta'n} \log   \mathbb{E}_n^V \left[ \exp \left(n\beta'\int (L\circ f)\, d\mu_V  + n\beta' \varphi( \zeta^{(n)},\gamma^{(n)} )\right)\right] \notag \\
& = \int (L\circ f)\, \diff \mu_V  + J(\varphi) .
\end{align}
The second equality follows from Theorem \ref{thm:LDPprojectivelimit} and Varadhans Lemma (Theorem 4.3.1 in \cite{demboz98}), with
\begin{align*}
J(\varphi) = \sup_{(z,g) \in \mathcal{Z}\times \mathcal{G} } \{ \varphi(z,g) - \mathcal{I}^{(\mathrm{ext},\mathrm{w})}(z,g) \}.
\end{align*}
Note that then by duality, also $\mathcal{I}^{(\mathrm{ext},\mathrm{w})}(z,g) = \sup_{\varphi \in C_b(\mathcal{Z}\times \mathcal{G})} \{ \varphi(z,g) - J(\varphi) \} $. This shows that the first assumption in Theorem \ref{newgeneral} holds, with $\Lambda(f) = \int L\circ f\, d\mu_V$. It was shown in \cite{magicrules} that 
\begin{align*}
\Lambda^*(\mu) = \mathcal{K}(\mu_V\!\ |\!\ \mu) +\mu(I) - 1 
\end{align*}
for $\mu \in \mathcal{M}(I)$. Moreover, the set $\mathcal{F}$ of exposed points of $\Lambda^*$ contains the set of measures $\mu=h\cdot \mu_V$, absolutely continuous with respect to $\mu_V$ with strictly positive continuous density $h$. The exposing hyperplane of  $\mu=h\cdot \mu_V$ is given by $1-h^{-1}$, such that for any such $\mu$ there exists a $\gamma>1$ such that $\gamma (1-h^{-1})\in D$. 
Suppose now that $\mu \in \mathcal{M}(I)$ is such that $\Lambda^*(\mu)$ is finite. 
By the same arguments as in \cite{grz}, we can find a sequence $\mu_n$ of measures with strictly positive continuous density such that $\mu_n$ converges weakly to $\mu$ and $\Lambda^*(\mu_n)$ converges to $\Lambda^*(\mu)$. This approximation is also made more precise for matrix valued measures in \cite{GaNaRomat}. 
All assumptions of Theorem \ref{newgeneral} are then fulfilled, which yields the joint LDP for $ (\tilde \mu_{n,I}, \zeta^{(n)},\gamma^{(n)})_n$. 
\hfill $\Box $

\subsection{Normalizing and recovering the spectral measure}
\label{sec:normalizing}

To finish the proof of Theorem \ref{MAIN}, two steps remain. First, we need to map the collection $(\tilde\mu_{n,I},\zeta^{(n)}, \gamma^{(n)})$ to the measure $\tilde\mu_n$, and then normalize $\tilde\mu_n$ to recover the distribution of the original spectral measure $\mu_n$. 

For the first step, let 
$\Theta:\mathcal{M}(I)\times \mathcal{Z}\times \mathcal{G}\rightarrow \mathcal{S}(I)$ be defined by
\begin{align} \label{defTheta}
\Theta(\mu_I,\zeta,\gamma) = \mu_I + \sum_{i=1}^{2M}\sum_{j=1}^\infty \gamma_{i,j} \delta_{\zeta_{i,j}} .
\end{align}
Then by the construction in Section \ref{sec:newencoding}, in particular \eqref{newencoding4}, we have 
\begin{align*}
\Theta(\tilde\mu_{n,I},\zeta^{(n)}, \gamma^{(n)}) = \tilde\mu_n .
\end{align*}
However, we cannot apply the contraction principle directly, since the mapping $\Theta $ is not continuous when $\mathcal{G}$ is endowed with product topology. 
 We need to slightly modify the LDP for $(\zeta^{(n)})_n$. Since the rate function for $(\zeta^{(n)})_n$ is given by the $\ell_1$-norm of an array in $\mathcal{G}$, it is easy to see that $(\zeta^{(n)})_n$ is exponentially tight in the $\ell_1$-topology. From \cite{demboz98}, Corollary 4.2.6 (and the LDP in Theorem \ref{thm:jointjointLDP}), we get that $(\tilde\mu_{n,I},\zeta^{(n)},\gamma^{(n)})_n$ satisfies under $\Pnv$ the LDP in $\mathcal{M}(I)\times \mathcal{Z}\times \mathcal{G}$, with $\mathcal{G}$ endowed with the $\ell_1$-topology, with speed $\beta'n$ and good rate function
$\widetilde{\mathcal{I}}$. We can then make use of the following lemma, the proof is postponed to the end of this section. 

\medskip

\begin{lem} \label{lem:continuous}
When $\mathcal{M}(I)$ is endowed with the weak topology, $\mathcal{Z}$ with the product topology, and $\mathcal{G}$ with the $\ell_1$-topology, the mapping $\Theta$ as defined in \eqref{defTheta} is continuous. 
\end{lem}

\medskip

Then by the contraction principle, the spectral measures $\tilde\mu_n = \Theta(\tilde\mu_{n,I},\zeta^{(n)}, \gamma^{(n)})$ satisfy under $\Pnv$ the LDP in $\mathcal{S}(I)$ with speed $\beta' n $ and good rate function 
\begin{align}\label{rateinS}
\widetilde{\mathcal{I}}_{\operatorname{sp}}(\tilde \mu) = \inf \left\{ \widetilde{\mathcal{I}}(\tilde{\mu}_I,z,g) \mid \Theta(\tilde{\mu}_I,z,g)=\tilde\mu \right\} . 
\end{align} 
Note that $\Theta$ is not a bijection: if $\tilde{\mu}$ has point masses in $\partial I$, they may come from $\tilde{\mu}_I$ or from elements of $g$, for which the corresponding entry in $z$ lies in $\partial I$, and a point mass of $\tilde\mu$ at $x\notin I$ may arise from the combination of several equal elements in $g$. It follows from the form of the rate $\widetilde{\mathcal{I}}$, that in the first case the infimum in $\eqref{rateinS}$ is obtained by attributing these point masses to $\tilde\mu_I$ and in the second case the infimum is attained by choosing only a single outlier at $x$. The infimum in \eqref{rateinS} is therefore given by 
\begin{align} \label{rateinS2}
\widetilde{\mathcal{I}}_{\operatorname{sp}}(\tilde \mu) & = \mathcal{K}(\mu_V | {\tilde\mu}) + \tilde\mu(I)-1+  \sum_{z \in E(\tilde\mu)} \mathcal{F}_V(z) + \tilde\mu(\{z \}) \notag \\
& = \mathcal{K}(\mu_V | {\tilde\mu}) + \tilde\mu(\mathbb{R})-1+  \sum_{z \in E(\tilde\mu)} \mathcal{F}_V(z)  .  
\end{align}
It remains to normalize the measures $\tilde\mu_n$. Note that if $\tilde\mu$ is the zero measure, the Kullback-Leibler part in \eqref{rateinS2} equals $+\infty$ and so the rate $\widetilde{\mathcal{I}}_{\operatorname{sp}}$ can only be finite if $\tilde\mu(\mathbb{R})>0$. Furthermore, $\Pnv(\tilde\mu_n(\mathbb{R})>0)=1$. Then we may restrict the LDP for $(\tilde\mu_n)_n$ to the set of measures $\tilde\mu \in \mathcal{S}(I)$ with $\tilde\mu(\mathbb{R})>0$ (see Lemma 4.1.5 in \cite{demboz98}). On this set of measures, the mapping 
$  \tilde\mu \mapsto\tilde\mu(\mathbb{R})^{-1}\tilde\mu$ 
%$\tilde\mu\mapsto \tilde\mu \cdot \tilde\mu(\mathbb{R})^{-1}$
 is continuous. Since $\tilde\mu_n(\mathbb{R})^{-1}\tilde\mu_n$ has the same distribution as 
%the distribution of $ \tilde\mu_n \cdot \tilde\mu_n(\mathbb{R})^{-1} $ equals the distribution of
 $\mu_n$, a final application of the contraction principle yields that $(\mu_n)_n$ satisfies the LDP in $\mathcal{S}_1(I)$ with speed $\beta'n$ and good rate function 
\begin{align} \label{finalrate}
\mathcal{I}_{\operatorname{sp}}(\mu) & = \inf_{\kappa>0} \widetilde{\mathcal{I}}_{\operatorname{sp}}(\kappa\mu%\mu\cdot \kappa
) \notag \\
& = \inf_{\kappa>0} \int \log \left( \frac{\mathrm{d} \mu_V}{\mathrm{d}(\kappa\mu%\mu \cdot \kappa 
)} \right) \mathrm{d}\mu_V
 + (\kappa\mu%\cdot \kappa
)(\mathbb{R})-1 + \sum_{z \in E(\tilde\mu)} \mathcal{F}_V(z) \notag \\
 & = \inf_{\kappa>0} (\kappa-\log \kappa -1 )+ \int \log\left( \frac{\mathrm{d} \mu_V}{\mathrm{d}\mu} \right) \mathrm{d}\mu_V
 +  \sum_{z \in E(\tilde\mu)} \mathcal{F}_V(z) .
\end{align}
This last infimum equals 0, attained for $\kappa=1$. This yields precisely the rate function in Theorem \ref{MAIN}. 

Finally, we can extend the last LDP for $(\mu_n)_n$ from the space $S_1(I)$ to $\mathcal{M}_1(\mathbb{R})$ by setting $\mathcal{I}_{\operatorname{sp}}(\mu)=+\infty$ if $\mu \notin S_1(I)$. Then it is easy to see that $\mathcal{I}_{\operatorname{sp}}$ is lower semicontinuous on $\mathcal{M}_1(\mathbb{R})$ and so the LDP holds in $\mathcal{M}_1(\mathbb{R})$ as well. This concludes the proof of Theorem \ref{MAIN}. 
\hfill $\Box$

\medskip

\textbf{Proof of Lemma \ref{lem:continuous}:}  
Let $\mu_{I}^{(n)}\rightarrow \mu_I $ weakly in $\mathcal{M}(I)$, $z^{(n)}\to z$ entrywise in $\mathcal{Z}$ and $g^{(n)}\to g$ in $\mathcal{G}$ with respect to the $\ell_1$-topology. Denote $\Theta(\mu_{I}^{(n)},z^{(n)},g^{(n)} )=\mu^{(n)}$ and $\Theta(\mu_{I},z,g )=\mu$. Let $f$ be continuous and bounded. Then 
\begin{align*}
\left| \int f \diff \mu^{(n)} - \int f \diff \mu \right| 
& \leq \left| \int f \diff \mu^{(n)}_I - \int f \diff \mu_I \right| 
		+ \sum_{i=1}^{2M} \sum_{j=1}^\infty | g^{(n)}_{i,j} f(z^{(n)}_{i,j}) - g_{i,j} f(z_{i,j})|  \\
& \leq \left| \int f \diff \mu^{(n)}_I - \int f \diff \mu_I \right| 
		+ \sum_{i=1}^{2M} \sum_{j=1}^q | g^{(n)}_{i,j} f(z^{(n)}_{i,j}) - g^{(n)}_{i,j} f(z_{i,j})|  \\
& \qquad + ||f||_\infty \sum_{i=1}^{2M} \sum_{j=q+1}^\infty |g_{i,j}| + ||f||_\infty ||g^{(n)}-g||_1  
\end{align*}
for any $q\geq 1$. The terms in the last two lines can be made arbitrarily small by first choosing $q$ and then $n$ large enough. 
\hfill $\Box$

\section{Proof of the coefficient LDP}
\label{sec:coeffproof}

The proofs of Theorem \ref{MAIN2} and Theorem \ref{MAIN3} make use of the explicit density \eqref{krishnadensity}, but for several arguments we rely on the fact that by Theorem \ref{MAIN}, we know an LDP holds for the spectral measure. In Section  \ref{susec:conditionalLDP}, this allows to show that an LDP holds for the recursion coefficients when we condition on a compact set $\mathcal{R}_K$. Although general large deviation theory allows to write the corresponding rate function as a projective limit, at this stage, it is not available in an explicit form. In Section \ref{susec:alternativerate}, we look at the density \eqref{krishnadensity} to obtain an alternative description for the rate function, up to an {\it error term}, which is bounded on the compact set $\mathcal{R}_K$. This proves Theorem \ref{MAIN2}. Finally, in Section \ref{susec:proofofonecut}, we show that in the one-cut case the error term vanishes, concluding the proof of Theorem \ref{MAIN3}.

\subsection{An abstract LDP for the conditional measure}
\label{susec:conditionalLDP}

To start with, note that by Theorem \ref{MAIN}, the sequence $(\mu_n)_n$ satisfies under $\mathbb{P}_n^V$ the LDP in $\mathcal{M}_1$ with speed $n\beta'$ and good rate function $\mathcal{I}_{\operatorname{sp}}$ which vanishes only at the compactly supported equilibrium measure $\mu_V$. The following theorem shows that this LDP holds also under conditioning on the smaller set $\mathcal{M}_{1,K}$ of probability measures with support in $[-K,K]$, where $K$ is so large that $I\subset [-K+1,K-1]$. We denote by $\mathbb{P}_{n,K}^V = \mathbb{P}_{n}^V(\cdot | \mu_n \in \mathcal{M}_{1,K})$ the measure conditioned on $\mathcal{M}_{1,K}$. 

\begin{thm} \label{thm:restrictedLDP}
Assume that the potential $V$ satisfies the assumptions (A1), (A2) and (A3). Then the sequence of spectral measures $\mu_n$ satisfies under $\mathbb{P}_{n,K}^V$ the LDP in $\mathcal{M}_{1,K}$ with speed $\beta'n$ and good rate function
the restriction of $\mathcal{I}_{\operatorname{sp}}$ to $\mathcal{M}_{1,K}$. 
\end{thm}

\proof
Instead of starting from Theorem \ref{MAIN}, we make use of Theorem \ref{thm:jointjointLDP}, which states that $(\tilde\mu_{n,I},\zeta_n,\gamma_n)_n$ satisfies the LDP in $\mathcal{M}(I)\times \mathcal{Z}\times \mathcal{G}$ with speed $\beta'n$ and good rate function $\mathcal{I}_{\operatorname{sp}}$. 
Let $\mathcal{Z}_K = \{ z \in \mathcal{Z}: \, ||z||_\infty \leq K\}$. Then $\mathcal{Z}_K$ is a closed subset of $\mathcal{Z}$. 
%(note that restricting the support in $\mathcal{M}(\mathbb{R})$ would not yield a closed subset). 
Furthermore, by the LDP in Theorem \ref{thm:extremalLDP} for the extremal eigenvalues, we have that $\Pnv(\zeta_n \in \mathcal{Z}_K)$ converges to 1. Therefore, for any set $C$ closed in $\mathcal{M}(I)\times \mathcal{Z}_K\times \mathcal{G}$, 
\begin{align} \label{restrictedLDPub}
\limsup_{n\to \infty} \frac{1}{\beta'n} \log \mathbb{P}_{n,K}^V ((\tilde\mu_{n,I},\zeta_n,\gamma_n) \in C) & = 
\limsup_{n\to \infty} \frac{1}{\beta'n} \log \mathbb{P}_{n}^V ((\tilde\mu_{n,I},\zeta_n,\gamma_n) \in C, ||\zeta_n||_\infty \leq K ) \notag \\
& \leq - \inf_{(\tilde\mu,z,g)\in C\cap \mathcal{M}(I)\times \mathcal{Z}_K\times \mathcal{G} } \widetilde{\mathcal{I}}(\tilde\mu,z,g) ,
\end{align} 
by the large deviation upper bound of Theorem \ref{thm:jointjointLDP}. 
Similarly, we get from the lower bound for any set $O$ open in $\mathcal{M}(I)\times \mathcal{Z}_K\times \mathcal{G}$, 
\begin{align} \label{restrictedLDPlb}
& \liminf_{n\to \infty} \frac{1}{\beta'n} \log \mathbb{P}_{n,K}^V ((\tilde\mu_{n,I},\zeta_n,\gamma_n) \in O) \notag  \\
 & \geq
\liminf_{n\to \infty} \frac{1}{\beta'n} \log \mathbb{P}_{n}^V ((\tilde\mu_{n,I},\zeta_n,\gamma_n) \in O\cap \mathcal{M}(I)\times \operatorname{Int}(\mathcal{Z}_K)\times \mathcal{G}) \notag \\
& \geq  - \inf_{(\tilde\mu,z,g)\in O\cap \mathcal{M}(I)\times \operatorname{Int}(\mathcal{Z}_K)\times \mathcal{G}} \widetilde{\mathcal{I}}(\tilde\mu,z,g) . 
\end{align}
where $\operatorname{Int}(\mathcal{Z}_K)$ is the interior of $\mathcal{Z}_K$ as a subset of $\mathcal{Z}$, that is, $\operatorname{Int}(\mathcal{Z}_K)= \{g \in \mathcal{Z}:\, ||g||_\infty <K\}$. We remark that this argument would not be helpful is we started from the LDP in $\mathcal{M}_1$, as then the interior of the restricted space (in the weak topology) would be empty. 

From the explicit form of the rate in Theorem \ref{thm:extremalLDP}, it can be seen that for any open set $O$, 
\begin{align*}
\inf_{(\tilde\mu,z,g)\in O\cap \mathcal{M}(I)\times \operatorname{Int}(\mathcal{Z}_K)\times \mathcal{G}} \widetilde{\mathcal{I}}(\tilde\mu,z,g)
= \inf_{(\tilde\mu,z,g)\in O\cap \mathcal{M}(I)\times \mathcal{Z}_K\times \mathcal{G}} \widetilde{\mathcal{I}}(\tilde\mu,z,g) . 
\end{align*}
Together with \eqref{restrictedLDPub}, this shows that $(\tilde\mu_{n,I},\zeta_n,\gamma_n)_n$ satisfies under $\mathbb{P}_{n,K}^V$ the LDP in the space $\mathcal{M}(I)\times \mathcal{Z}_K\times \mathcal{G}$ with rate function the restriction of $\widetilde{\mathcal{I}}$. 

We may now proceed as in the proof of Theorem \ref{MAIN}. We have $\mu_n \in \mathcal{M}_{1,K}$ if and only if $\zeta_n \in \mathcal{Z}_K$. The same arguments as in  Section \ref{sec:normalizing} applied to the restricted LDP show that $(\mu_n)_n$ satisfies the LDP in the space $\mathcal{M}_{1,K}$, and the rate function is the restriction of $\mathcal{I}_{\operatorname{sp}}$ to this space. 
\hfill $\Box$

\begin{cor} \label{cor:restrictedLDP}
Assume that the potential $V$ satisfies the assumptions (A1), (A2) and (A3). Then the sequence of recursion coefficients $r_n$ satisfies under $\mathbb{P}_{n,K}^V$ the LDP in $\mathcal{R}_K$ with speed $\beta'n$ and good rate function given by 
\begin{align*}
\mathcal{I}_{\operatorname{co},K}(r) = \mathcal{I}_{\operatorname{sp}}(\psi^{-1}(r)) = \lim_{N\to \infty} \mathcal{I}_N(r),  
\end{align*}
with
\begin{align*}
\mathcal{I}_{N}(r) = -  \lim_{\delta \to 0} \limsup_{n\to \infty} \frac{1}{\beta'n} \log \mathbb{P}_{n,K}^V (B_{\delta,N} (r)) , 
\end{align*}
where $B_{\delta,N} (r)= \{ z \in \mathbb{R}^\mathbb{N} |\, |z_i-r_i|<\delta \text{ for } i\leq 2N-1\}$ is the ball around the first $2N-1$ coordinates of $r$.  
\end{cor}

\proof 
We have $\mathcal{R}_K = \psi(\mathcal{M}_{1,K})$, and $\psi$ is a homeomorphism from $\mathcal{M}_{1,K}$ to $\mathcal{R}_K$, which implies by the contraction principle the LDP for $(r_n)_n$ with good rate function $\mathcal{I}_{\operatorname{sp}}\circ \psi^{-1}$. Restricting the continuous projections to $\mathcal{R}_K$, we get again by the contraction principle that the sequence of projected coefficients $\pi_N(r_n)$ satisfies the LDP in $\pi_N(\mathcal R_K)$, with some good rate function $\tilde{\mathcal{I}}_N$. The Dawson-G\"artner Theorem implies that the rate function for $(r_n)_n$ can then be recovered as  
$\mathcal{I}_{\operatorname{co},K}= \lim_{N\to \infty} \tilde{\mathcal{I}}_N\circ \pi_N$. On $\mathcal R_K$, we let $\mathcal{I}_N=\tilde{\mathcal{I}}_N\circ \pi_N$. As shown in Theorem 4.1.18 in \cite{demboz98}, 
\begin{align}
\tilde{\mathcal{I}}_N(\pi_N(r)) = -  \lim_{\delta \to 0} \limsup_{n\to \infty} \frac{1}{\beta'n}  \mathbb{P}_{n,K}^V (||\pi_N(r_n)- \pi_N(r)||_\infty < \delta ) , 
\end{align}
which proves the last display of Corollary \ref{cor:restrictedLDP}. \hfill $ \Box$

\medskip

In the following, we write $r_{n,N}$ for $\pi_N(r_n)$, and if $T_n$ is the tridiagonal matrix with Jacobi coefficients $r_n$, we write $T_{n,N}$ for $\pi_N(T_n)$. 
 Recall that $r^V$ is the sequence of Jacobi coefficients of $\mu_V$ and the corresponding Jacobi operator is $T^V$. We use the analogous notation for $r^V_N=\pi_N(r^V)$ and $T^V_N=\pi_N(T^V)$.

\subsection{An alternative expression for the rate}
\label{susec:alternativerate}

To obtain an alternative description of $\mathcal{I}_N$, we will decompose the density in \eqref{krishnadensity} into three factors, one depending only on $r_{n,N}$, one only on entries omitted in $r_{n,N}$ and one factor containing finitely many mixed terms. Let $d=2p$ be the degree of the polynomial potential $V \in \mathcal V$.

\begin{lem}\label{lem:densitydecomp}
There exist continuous functions $M : ([0, \infty) \times \R)^{2d+1} \to \R$ not depending on $n$ and $E_n : ([0, \infty) \times \R)^{n-N} \to \R$, such that for all $n\geq N+ d, N\geq d$,
\begin{align*}
\tr\!\ V(T_n)-\tr\!\ V(T^V_n) = \ & \tr\!\ V(T_{n,N})- \tr\!\ V(T^V_{N}) + M(a_{N-d},b_{N-d},\dots ,b_{N+d})\\ 
&+ E_n(a_{N+1},b_{N+1},\dots,b_n )\,.
\end{align*}
Moreover, if $|a_k|,|b_k|\leq K$ for every $k \leq n$, then there exists a constant $C(K,V) > 0$ such that for every $N\geq d$ : 
\begin{align*}
|M(a_{N-d},b_{N-d},\dots , b_{N+d})| \leq  C(K,V) M_+(r_N)\,,
%\big( |a_{N-d}-a^V_{N-d}| + |b_{N-d}-b^V_{N-d}| + \dots +  |b_{N}-b^V_{N}|\big) .
\end{align*}
with $M_+$ defined as in \eqref{defM+}.
\end{lem}

\textbf{Proof:} By linearity, it suffices to show the decomposition for $V(x)=x^d$ a monomial. Note that $\tr V(T_{n,N})= \tr V(A)$, where $A = T_{n,N} \oplus 0_{n-N}$ and $T_{n,N}$ is the $N\times N$ tridiagonal matrix with the first $2N-1$ entries of $r_n$. Let $B= T_n- A$. We have
\begin{align*}
T_n^d = (A+B)^d = A^d + B^d + \sum_{i\in \{0,1\}^d, i\neq 0,1} A^{i_1}B^{1-i_1} \dots B^{1-i_d} ,
\end{align*}
where in the last sum there is always one factor equal to $A$ and one equal to $B$. Define $\hat A$, $\hat B$ analogously, build from $T_n^V$, then
\begin{align*}
V(T_n)-V(T_n^V) & = (A^d-\hat A^d) + (B^d- \hat B^d) +   \sum_{i\in \{0,1\}^d, i\neq 0,1} \left( A^{i_1}B^{1-i_1} \dots B^{1-i_d} - \hat A^{i_1}\hat B^{1-i_1} \dots \hat B^{1-i_d} \right) \\
& = (A^d-\hat A^d) + (B^d+\hat B^d) +   \sum_{i\in \{0,1\}^d, i\neq 0,1} \left( A^{i_1}B^{1-i_1} \dots B^{1-i_d}- \hat A^{i_1}B^{1-i_1} \dots B^{1-i_d}\right) \\
& \qquad  +  \sum_{i\in \{0,1\}^d, i\neq 0,1} \left( \hat A^{i_1}B^{1-i_1} \dots B^{1-i_d} - \hat A^{i_1}\hat B^{1-i_1} \dots \hat B^{1-i_d} \right) .
\end{align*}
Now $(A^d-\hat A^d)=  V(T_{n,N})- V(T^V_{n,N})$, and on the other hand $(B^d-\hat B^d)$ and the last sum do not depend on $r_{n,N}$ and their trace can be combined into $E_n$. We are then left with evaluating 
\begin{align*}
%\tilde M 
\Delta=   \sum_{i\in \{0,1\}^d, i\neq 0,1} \left( A^{i_1}B^{1-i_1} \dots B^{1-i_d}- \hat A^{i_1}B^{1-i_1} \dots B^{1-i_d}\right) .
\end{align*}
Suppose $n\geq d$. To see that $\tr \Delta$ %\tilde M$
 depends only on $a_{N-d},b_{N-d},\dots ,b_{N+d}$, write
\begin{align} \label{eq:traceissues}
\tr \left( A^{i_1}B^{1-i_1} \dots B^{1-i_d}\right) =
 \sum_{k_1,\dots ,k_d=1}^n   (A^{i_1}B^{1-i_1})_{k_1,k_2}(A^{i_2}B^{1-i_2})_{k_2,k_3} \dots(A^{i_d} B^{1-i_d})_{k_d,k_1} .
\end{align}
Both $A$ and $B$ are tridiagonal, such that any nonzero term in this sum satisfies $|k_i-k_{i-1}|,|k_1-k_d|\leq 1$. In other words, $k=(k_1,\dots, k_d,k_1)$ is a closed path on $\{1,\dots ,n\}$ with step size at most 1. Furthermore, $A_{\ell,m}=0$ if $\ell\geq N+1$ or $m\geq N+1$ and $B_{\ell,m}=0$ if $\ell\leq N-1$ or $m\leq N-1$. At least one of the matrices $(A^{i_j}B^{1-i_j})$ equals $A$ and one equals $B$. Therefore, any path $k=(k_1,\dots, k_d,k_1)$ with $k_i\neq N$ for all $i$ gives a zero term in \eqref{eq:traceissues}. But then any contribution in \eqref{eq:traceissues} comes from paths with $N-\lfloor d/2\rfloor \leq k_1,\dots ,k_d\leq N+\lfloor d/2\rfloor$, which implies that only the entries $a_{N-d},b_{N-d},\dots ,b_{N+d}$ appear in \eqref{eq:traceissues}. The same holds true if we replace $A$ by $\hat A$, so that $M=\tr\Delta$ has the claimed form. 

It remains to show the bound for $|M|$. After taking the trace, we are left with finitely many differences 
\begin{align*}
& (A^{i_1}B^{1-i_1})_{k_1,k_2}(A^{i_2}B^{1-i_2})_{k_2,k_3} \dots(A^{i_d} B^{1-i_d})_{k_d,k_1} \\
& \qquad \qquad - (\hat A^{i_1}B^{1-i_1})_{k_1,k_2}(\hat A^{i_2}B^{1-i_2})_{k_2,k_3} \dots(\hat A^{i_d} B^{1-i_d})_{k_d,k_1} ,
\end{align*}
with $k$ a path as above. Whenever one of the entries of $r_{n,N}$ appears in the first product (and there is always at least one such entry), the corresponding entry of $r^V_{n,N}$ appears in the second product, and the desired bound follows from the boundedness of $|a_k|,|b_k|$ and possibly the triangle inequality, in case $A$ appears more than once in the product. 
\hfill $\Box$ 
\medskip

Looking at the density \eqref{krishnadensity}, a natural guess for the rate function of the projected vector $r_{n,N}$ would be
\begin{align}
\label{defUJ}
\mathcal{U}_N(r_N) = \tr\!\ V(T_N) -\tr\!\ V(T^V_N) - 2\sum_{k=1}^{N-1} \log (a_k/a^V_k) .
\end{align}
Since we cannot ignore the boundary effects with higher order Jacobi coefficients in \eqref{krishnadensity}, we cannot conclude the LDP with this rate function. In fact, the deviation from $\mathcal{U}_N$ will be given in terms of $M_+(r_N)$.
%the bound for the function $M$ in Lemma \ref{lem:densitydecomp}, for which we define
%\begin{align*}
%M_+(r_J) =  C(K,V)\big( |a_{J-d}-a^V_{J-d}| + |b_{J-d}-b^V_{J-d}| + \dots +  |b_{J}-b^V_{J}|\big) .
%\end{align*}
We then have the following result.

\begin{thm} \label{thm:improperLDP}
Let $r_N$ be a fixed finite vector in $\pi_N(\mathcal{R}_K)$ and let $B_\delta(r_N)$ be the open ball in $\mathbb{R}^{2N-1}$ around $r_N$ with respect to the sup-norm. Then 
\begin{align*}
\lim_{\delta \to 0} \limsup_{n\to \infty} \frac{1}{n\beta'} \log \mathbb{P}_{n,K}(r_{n,N} \in B_\delta(r_N)) & \leq - \mathcal{U}_N(r_N) + C(K,V) M_+(r_N) \\
\lim_{\delta \to 0} \liminf_{n\to \infty} \frac{1}{n\beta'} \log \mathbb{P}_{n,K}(r_{n,N} \in B_\delta(r_N)) & \geq - \mathcal{U}_N(r_N) - C(K,V) M_+(r_N) \\
\end{align*}
\end{thm}

\textbf{Proof:} 
We use the same idea as in \cite{BSZ} and look at ratios of probabilities, so that we can ignore the normalizing constant and consider $\widetilde{\mathbb{P}}^V_{n,K}=Z^V_{n,K}\mathbb{P}^V_{n,K}$. We then decompose the density as in Lemma \ref{lem:densitydecomp}. For this, define additionally
\begin{align*}
\ell_0(a_1,\dots ,a_{N-1}) & = 2 \sum_{k=1}^{N-1} \left( \frac{k}{n} + \frac{1}{n\beta}\right) \log (a_k/a^V_k) , \\
\ell_1(a_N,\dots ,a_{n-1}) & = -2 \sum_{k=N}^{n-1}  \left( 1- \frac{k}{n} - \frac{1}{n\beta}\right) \log (a_k/a^V_k) .
\end{align*}
The density of Jacobi coefficients is given in \eqref{krishnadensity}. The measure $\widetilde{\mathbb{P}}^V_{n,K}$ has a (non-normalized) density, which on $\pi_N(\mathcal{R}_K)$ is proportional to \eqref{krishnadensity}. The restriction to $\pi_N(\mathcal{R}_K)$ implies in particular, that the density of $\widetilde{\mathbb{P}}_{n,K}$ is zero on the complement of $[-K,K]\times ([0,K]\times [-K,K])^{n-1}$. 
Given the ball $B_\delta(r_N)$, define 
\begin{align*}
\widetilde B_\delta(r_N) = B_\delta(r_N) \times ([0,K]\times [-K,K])^{n-N} .
\end{align*}
 Then, using the decomposition from Lemma \ref{lem:densitydecomp},
\begin{align*}
\widetilde{\mathbb{P}}^V_{n,K} (\widetilde B_\delta(r_N) ) = \int_{\widetilde B_\delta(r_N)\cap \mathcal{R}_K} \exp\left\{-n \beta'\left( \mathcal{U}_N + M+E_n+\ell_0+\ell_1\right) \right\} d\lambda_n ,
\end{align*}
with $E_n$ and $\ell_1$ independent of $b_1,a_1,\dots b_j$, and $\mathcal{U}_N$ and $\ell_0$ independent of $a_{N},b_{N+1},\dots ,b_n$. Here, we wrote $\lambda_n$ for the Lebesgue measure on $\mathbb{R}^{2n-1}$. 
Looking at the ratio of probabilities and applying the bound for $M$ in Lemma \ref{lem:densitydecomp}, we have then 
\begin{align*}
& \quad \frac{1}{n\beta'} \log \frac{{\mathbb{P}}^V_{n,K} (\widetilde B_\delta(r_N) )}{{\mathbb{P}}_{n,K}^V (\widetilde B_\delta( r^V_N) )} 
 = \frac{1}{n\beta'} \log \frac{\widetilde{\mathbb{P}}^V_{n,K} (\widetilde B_\delta(r_N) )}{\widetilde{\mathbb{P}}^V_{n,K} (\widetilde B_\delta(r^V_N) )} \\
& \leq \sup_{ r\in B_\delta(r_N)}  \left(  - \mathcal{U}_N(r) - \ell_0(r) + C(K,V) M_+(r)\right) - \inf_{ r\in B_\delta( r^V_N)}  \left( -  \mathcal{U}_N(r) - \ell_0(r) - C(K,V) M_+(r)\right) .
\end{align*}
By continuity of $\mathcal{U}_N,\ell_0$ and $M_+$ on $B_\delta(r_N)$, 
\begin{align*}
\lim_{\delta \to 0} \lim_{n\to \infty} \sup_{ r\in B_\delta(r_N)}  \left(  - \mathcal{U}_N(r) - \ell_0(r) +C(K,V) M_+(r)\right)  
&= \lim_{\delta \to 0} \sup_{ r\in B_\delta(r_N)}\left( -\mathcal{U}_N(r) + C(K,V)M_+(r)\right) \\
& = - \mathcal{N}_N(r_N) +C(K,V) M_+(r_N),
\end{align*}
and 
\begin{align*}
\lim_{\delta \to 0} \lim_{n\to \infty} \inf_{ r\in B_\delta(r_N)}  \left(  - \mathcal{U}_N(r) - \ell_0(r) - C(K,V) M_+(r)\right) 
& = \lim_{\delta \to 0} \sup_{ r\in B_\delta(r_N)}\left( -\mathcal{U}_N(r) - C(K,V) M_+(r)\right) \\
& = -\mathcal{U}_N(r_N) - C(K,V) M_+(r_N).
\end{align*}
For the ratio of probabilities this implies
\begin{align*}
\lim_{\delta\to 0} \limsup_{n\to \infty} \frac{1}{n\beta'} \log \frac{{\mathbb{P}}^V_{n,K} (\widetilde B_\delta(r_N) )}{{\mathbb{P}}^V_{n,K} (\widetilde B_\delta( r^V_N) )} 
& \leq -\mathcal{U}_N(r_N) +C(K,V) M_+(r_N)+\mathcal{U}_N(r^V_N) +C(K,V) M_+( r^V_N) \\
& = -\mathcal{U}_N(r_N) +C(K,V) M_+(r_N)
\end{align*}
Since ${\mathbb{P}}^V_{n} (\widetilde B_\delta( r^V_N) )$ and then also ${\mathbb{P}}^V_{n,K} (\widetilde B_\delta( r^V_N) )$ converges to 1, this implies the first inequality, the second one follows by analogous arguments. 
\hfill $\Box$ 

\medskip

Theorem \ref{thm:improperLDP} implies for the rate function $\mathcal{I}_{N}$ of the projected sequence $\pi_N(r_n)$ 
\begin{align} \label{comparison}
\big| \mathcal{I}_{N}(r_N) -  \mathcal{U}_N(r_N) \big| \leq  C(K,V) M_+(r_N) 
\end{align}
for any $r_N \in \pi_N(R_K)$. By the second identity in Corollary \ref{cor:restrictedLDP}, the sequence $(r_n)_n$ satisfies then the LDP with speed $\beta'n$ and rate function given by $\mathcal{I}_{\operatorname{co}}(r)=\lim_{N\to \infty} \mathcal{I}_{N}(\pi_N(r))$. We may then set $\xi_{N,K} =  \mathcal{U}_N - \mathcal{I}_N$ and get for the rate of $r_n$ 
\begin{align} \label{projrate}
\mathcal{I}_{\operatorname{co}}(r)=\lim_{N\to \infty} \left[ \mathcal{U}_N(\psi_N(r)) + \xi_{N,K} \right]. 
\end{align}
Together with \eqref{comparison} and the bound %for $M_+$ 
in Lemma \ref{lem:densitydecomp}, this implies Theorem \ref{MAIN2}.

\subsection{Reduction to the one-cut case: proof of Theorem \ref{MAIN3}} 
\label{susec:proofofonecut}

First, suppose $r_n$ is distributed according to $\mathbb{P}^V_{n,K}$. 
We will show that the limit in \eqref{projrate} equals $\lim_{N\to \infty} \mathcal{U}_N(\pi_N(r))$, using the large deviation result of Theorem \ref{MAIN2}. Suppose that $\mathcal{I}_{\operatorname{co}}(r)$ is not finite for some $r\in \mathcal{R}_K$. Then \eqref{comparison} and the uniform bound for $M_+(\pi_N(r))$ on $\mathcal{R}_K$ implies that $\lim_{N\to \infty} \mathcal{U}_N(\pi_N(r))$ is infinite as well.  
Suppose $\mathcal{I}_{\operatorname{co}}(r)$ is finite. Since $r\in \mathcal{R}_K$ there exists a unique $\mu$ with support in $[-K,K]$ such that $\psi(\mu) = r$. By the contraction principle, $\mathcal{I}_{\operatorname{sp}}(\mu)<\infty$. By the Kullback-Leibler part of the rate, $\mu$ has then a Lebesgue decomposition $\diff \mu(x)= f(x) \diff\mu_V(x) + \diff\mu_s(x)$ with $f(x)>0$ for $\mu_V$-almost all $x\in \operatorname{supp}(\mu_V)$. By the explicit form of $\mu_V$ as in \eqref{eq:polynomialmuV}, this implies $f(x)>0$ for Lebesgue-almost all $x\in\operatorname{supp}(\mu_V)$. Rakhmanov's Theorem for Jacobi matrices \cite{Den2004} yields that then $a_k(\mu)\to \hat a$ and $b_k(\mu)\to \hat b$, where $\hat a = \lim_{k\to \infty} a_k(\mu_V), \hat b = \lim_{k\to \infty} b_k(\mu_V)$. From the bound for $M_+$ in Lemma \ref{lem:densitydecomp}, we have 
\begin{align} \label{Mvanishes}
\lim_{N\to \infty} M_+(\pi_N(r)) =0, 
\end{align}
and then 
\begin{align}\label{projlimit2}
\mathcal{I}(r) = \lim_{N\to \infty} \mathcal{U}_N(\pi_N(r)) 
\end{align} 
as well. 
It remains to extend the LDP to the full space $\mathcal{R}$ defined in \eqref{defR}. From \eqref{exptight2}, we have
\begin{align} \label{Qexptight3}
\lim_{K\to \infty} \limsup_{n\to \infty} \frac{1}{n\beta'} \log  \mathbb{P}^V_n(\mathcal{R}^c_K) = -\infty ,
\end{align}
such that the measures $\mathbb{P}^V_{n,K}$ are exponentially good approximations of the measures $\mathbb{P}^V_n$. By Theorem 4.2.16 in \cite{demboz98}, the sequence $(r_n)_n$ under $\mathbb{P}^V_n$ satisfies the LDP with speed $\beta'n$ and rate given by the limit of \eqref{projlimit2} as $K\to \infty$.

\section{LDP for extremal eigenvalues}
\label{sect:extremalLDP}

In this section we prove Theorem \ref{thm:extremalLDP}.
We first remark that $\mathcal I_N^{\ext}$ is a good rate function: it is lower semicontinuous as proved in \cite{BorGui2013onecut}, A.1. p.478.  From the same reference, 
$\Fr_{V}$ has compact level sets, so that  
$\mathcal I_N^{\ext}$ has compact level sets by the union bound. 
In Section \ref{sect:extremalLDP1}, we show exponential tightness of $(\zeta_N^{(n)})_{n\geq 1}$ under the sequence $\mathbb{P}_n^V$. It then suffices to prove the weak LDP, which follows from the control of probabilities of balls $B_\delta(z)$ of radius $\delta$ in the sup-norm around $z \in \mathcal{Z}_N$. We then show in Section \ref{sect:extremalLDP2} the upper bound 
\begin{align}
\label{ub}
\lim_{\delta \to 0} \limsup_{n\to \infty}\ (\beta'n)^{-1} \log \Pnv ( \zeta_N^{(n)} \in B_\delta(z) ) \leq -  \mathcal{I}_N^{\operatorname{ext}}(z)
\end{align}
and in Section \ref{sect:extremalLDP3} the lower bound 
\begin{align}
\label{lb}
\lim_{\delta \to 0} \liminf_{n\to \infty}\ (\beta'n)^{-1} \log \Pnv(\zeta_N^{(n)} \in B_\delta(z)) \geq - \mathcal{I}_N^{\operatorname{ext}}(z) ,
\end{align}
for any $z \in \mathcal{Z}_N$, which together imply then the full LDP. 
Along the way, we need the following four technical lemmas. Since their proofs are straightforward generalizations of the one-cut proof in \cite{magicrules}, they are omitted.

\begin{lem}
\label{LDPtilde}
Let  $V$ be  a potential satisfying the confinement condition (A1) and let $r$ be a fixed integer.  If  $\mathbb P^{V_n}_n$ is the probability measure associated to the potential
$V_n= \frac{n+r}{n} V$, 
then the law of $\mu_n^\u$ under $\mathbb P_n^{V_n}$ satisfies the LDP with speed $\beta' n^2$ with good rate function 
\begin{equation}  
\label{entropyV}
\mu \mapsto \mathcal E(\mu) - \inf_\nu\mathcal E(\nu)
\end{equation}
where $\mathcal E$  is defined in \eqref{ratemuu}.
\end{lem}

\begin{lem}
\label{lem:teknik}
If the potential $V$ is finite and continuous on a compact set and infinite outside, we have, for every $q \geq 1$
\begin{equation}
\lim_{n\to \infty} \frac{1}{n} \log \frac{Z^V_n}{Z^{\frac{n}{n-q}V}_{n-q}} = -q \inf_{x\in \mathbb{R}} \mathcal J_V(x)\,.
\end{equation}
\end{lem}

\begin{lem}\label{lem:cvproba}
Under Assumption (A1) and (A3), $\max_{i=1,\dots 2M} d(\zeta_{i,1}^{(n)} ,\partial I)$ converges to 0 in probability. Also, for any $q\geq 1$ and $\varepsilon>0$, 
\begin{align*}
\lim_{n\to \infty} \mathbb{P}_{n-q}^{\frac{n}{n-q}V} \left( \max_{i=1,\dots, 2M} d(\zeta_{i,1}^{(n)} ,\partial I) >\varepsilon \right) = 0 . 
\end{align*}
\end{lem}

\begin{lem} \label{lem:Zratio}
Under Assumption (A1),
\begin{align*}
\limsup_{n\to \infty} \frac{1}{n} \log \frac{Z_{n-1}^V}{Z_n^V}  < \infty . 
\end{align*}
\end{lem}

\subsection{Exponential tightness}
\label{sect:extremalLDP1}

The exponential tightness will follow from  
\begin{align}
\label{expotight}
\limsup_{L \rightarrow \infty}\limsup_{n \rightarrow \infty} \frac{1}{n} \log \Pnv\big( \zeta_N^{(n)} \notin K_M^{2MN} \big) = -\infty
\end{align}
for any $N\geq 1$, with $K_L = \{x\in \mathbb{R} \mid V(x)\leq L\}$. For $L$ large enough, we have
\begin{align} \label{expotight2}
\Pnv\big( \zeta_N^{(n)} \notin K_L^{2MN} \big)\leq \Pnv( \zeta_{1,1}^{(n)} \notin K_L) + \Pnv( \zeta_{2M,1}^{(n)} \notin K_L)
\end{align}
so the proof of exponential tightness reduces to the consideration of the smallest and the largest eigenvalue, and by symmetry, it suffices to show
\begin{align}
\label{expotight3}
\limsup_{L \rightarrow \infty}\limsup_{n \rightarrow \infty} \frac{1}{n} \log \Pnv(  \zeta_{2M,1}^{(n)} \notin K_L ) = -\infty . 
\end{align}
The rest of the proof follows now verbatim the proof of (A.7) in \cite{magicrules}, making use of Lemma \ref{lem:Zratio}.

As a consequence of exponential tightness, we may simplify the remaining proof substantially by replacing the potential $V$ by 
\begin{align*}
V_L (x) = \begin{cases} V(x) & \text{ if } V(x) \leq L ,\\
\infty & \text{ otherwise},
\end{cases}
\end{align*}
for $L$ large enough. Indeed, if $L$ is large enough, the minimizer $\mu_{V_L}$ will coincide with $\mu_V$ and also $\inf_{\xi\in \mathbb{R}} \mathcal J_{V_L}(\xi)=\inf_{\xi\in \mathbb{R}} \mathcal J_V(\xi)$. For the sake of a lighter notation, we will drop the subscript $L$, but we may assume that the eigenvalues are confined to a compact interval. In particular, Lemma \ref{lem:teknik} is applicable.

\subsection{Proof of the upper bound}
\label{sect:extremalLDP2}

In this section, we prove the upper bound \eqref{ub}.  Let $z\in \mathcal{Z}_N$. Without loss of generality, we may assume that $z_{i,j} \notin I$ for all $i,j$. To see this, let $\mathrm{Ind} = \{ (i,j): z_{i,j} \notin I\}$ be the set of indices of entries not in $I$. Then we have the trivial upper bound
\begin{align} \label{ubtrivial}
\Pnv ( \zeta_N^{(n)} \in B_\delta (z) ) \leq \Pnv ( \zeta_{i,j}^{(n)} \in [z_{i,j}-\delta,z_{i,j}+\delta] \text{ for all } (i,j)\in \mathrm{Ind} ) , 
\end{align}
and since $\mathcal{F}_V(z_{i,j}) = 0$ for $(i,j)\notin \mathrm{Ind}$, it suffices to consider the entries with indices in $\mathrm{Ind}$. In order to keep the notation simple, we assume then that $z_{i,j} \notin I$ for all $i,j$. In addition, let $\delta$ be so small that $B_\delta(z)\cap I^{2MN} = \emptyset$.

The eigenvalue density as in \eqref{generaldensity} is the density of unordered eigenvalues, so that we have  
\begin{align} \label{ub:representation0}
\Pnv (\zeta^{(n)}_N\in B_\delta(z)) 
 = \binom{n}{2MN} \frac{1}{Z^V_n} \int_{B_\delta(z)} \int_{\Delta(\lambda^{\mathrm{ex}})}    \prod_{i< j}|\lambda_i-\lambda_j|^\beta \prod_{i=1}^n e^{-\beta'V(\lambda_i)} \diff\lambda^\mathrm{in} \diff \lambda^\mathrm{ex} ,
\end{align}
where $\lambda^\mathrm{ex}=(\lambda_1,\dots ,\lambda_{2MN}) \in \mathbb{R}^{2MN}$ are the collection of (unordered) extremal eigenvalues, and the vector of the remaining (unordered) eigenvalues is denoted by $\lambda^\mathrm{in}=(\lambda_{2MN+1},\dots ,\lambda_n)\in \mathbb{R}^{n-2MN}$. Here, 
we consider $B_\delta(z)$ also as a subset of $\mathbb{R}^{2MN}$.  
Fixing the extremal eigenvalues forces the entries of $\lambda^\mathrm{in}$ then to be in the compact set 
\begin{align}
D(\lambda^\mathrm{ex}) = \left( \bigcup_{i=1}^{2M} \big[ \max \{\lambda_k : \, k\leq 2MN, \lambda_k\leq l_i\} , \min\{\lambda_k:\, k\leq 2MN, \lambda_k\geq r_i \} \big] \right)^{n-2MN} ,
\end{align}
where we recall that $I=[l_1,r_1]\cup \dots \cup [l_m,r_m]$. That is, the elements of $D(\lambda^\mathrm{ex})$ are ``more internal'' than the vector of eigenvalues $\lambda^\mathrm{ex}$, according to the ordering introduced in Section \ref{sec:newencoding}. For any $\lambda^\mathrm{ex}\in B_\delta(z)$, the maxima and minima in the definition of $\Delta(\lambda^{ex})$ are attained. 
The integral in \eqref{ub:representation0} may be rewritten as
\begin{align} \label{ub:representation}
\Pnv (\zeta^{(n)}_N\in B_\delta(z)) 
 = \binom{n}{2MN} \frac{1}{Z_n^V} \int_{B_\delta(z)} \Upsilon_{n,N}(\lambda^\mathrm{ex})\ d\lambda^\mathrm{ex}  ,
\end{align}
with the term $\Upsilon_{n,N}(\lambda^\mathrm{ex})$ given by 
\begin{align} \label{ub:upsilon} 
\Upsilon_{n,N}(\lambda^\mathrm{ex})=  H(\lambda^\mathrm{ex}) \Xi_{n,N}(\lambda^\mathrm{ex})  \exp\left\{ -\beta'n\sum_{k=1}^{2MN}V(\lambda_k) \right\} ,
\end{align}
with 
\begin{align*}
H(\lambda^\mathrm{ex}) =\prod_{1\leq r<s\leq 2MN} |\lambda_r - \lambda_s|^{\beta} ,
\end{align*}
and
\begin{align}  \label{ub:upsilon2}
\Xi_{n,N}(\lambda^\mathrm{ex}) & = \int_{D(\lambda^\mathrm{ex}) } \prod_{r=1}^{2MN}   \prod_{s=2MN+1}^{n} |\lambda_r-\lambda_s|^\beta  \prod_{r=2MN+1}^{n} e^{-n\beta'V(\lambda_r)}\prod_{2MN < r <s \leq n} |\lambda_r - \lambda_s|^\beta  \diff\lambda^\mathrm{in} \notag \\
& = Z_{n-2MN}^{\frac{n}{n-2MN}V}
\int_{D(\lambda^\mathrm{ex})} \prod_{r=1}^{2MN} \prod_{s=2MN+1}^{n} |\lambda_r-\lambda_s|^\beta   
 \diff \mathbb P_{n-2MN}^{\frac{n}{n-2MN}V} (\lambda) .
\end{align} 
In order to simplify notation we define $q=2MN$, so the above measure becomes $\mathbb P_{n-q}^{\frac{n}{n-q}V}$. 
%In the above integral, we may replace $V$ by
%\begin{align*}
%V_M (x) = \begin{cases} V(x) & \text{ if } x\in [-M,M] ,\\
%\infty & \text{ otherwise},
%\end{cases}
%\end{align*}
%for $M$ large enough, so that $V_M$ satisfies the assumption of Lemma \ref{lem:teknik}. Then we have
Now, to find an upper bound for $\Upsilon_{n,N}(\lambda^\mathrm{ex})$, we 
first choose   $K$ so large that $\operatorname{supp}(\mu_V)\subset [-K+1,K-1]$ and define 
\begin{align}\label{def:ball1}
\mathcal{B}_\kappa = \big\{ \mu \in \mathcal{M}_{1,K} |\, d_P(\mu,\mu_V)<\kappa \big\}
\end{align}
 the open ball around $\mu_V$ with radius $\kappa$ in the Prokhorov-metric, %assuming support in $[-K,K]$,
 %with $K$ so large that $\operatorname{supp}(\mu_V)\subset [-K+1,K-1]$ 
(recall the definition of $\mathcal M_{1,K}$ in (\ref{defMK})). Let also $\bar{\mathcal{B}}_\kappa=\{ \lambda \in \R^{n-q} |\,  \mu^\u_{n-q} \in \mathcal B_\kappa \}$. On the bounded set $D(\lambda^\mathrm{ex})$ the integrand in \eqref{ub:upsilon2} can be bounded by $ e^{c_1n}$ for some $c_1 > 0$ depending only on $z$ and $\delta$. 
We then use the fact that by Lemma \ref{LDPtilde}, the sequence of measures $\mu^\u_{n-q}$ satisfies under $\mathbb P_{n-q}^{\frac{n}{n-q}V}$ the LDP with speed $\beta'n^2$, and rate function vanishing only at $\mu_V$. The same arguments as in the large deviation upper bound in  \cite{magicrules} yield then
for any $\eta>0$,
\begin{align} \label{ub:ub}
& \limsup_{n\to \infty} \, (\beta'n)^{-1}\log \Pnv (\zeta_N^{(n)} \in B_\delta(z) ) \notag \\
&   = \limsup_{n\to \infty} \, (\beta'n)^{-1}\log \int_{B_\delta(z)} \left(Z^{\frac{n}{n-q}V}_{n-q}\right)^{-1} \Upsilon_{n,N}(\lambda^\mathrm{ex}) \diff \lambda^{\mathrm{ex}} 
				+ \limsup_{n\to \infty}\,  (\beta'n)^{-1} \log \left(\frac{Z^{\frac{n}{n-q}V}_{n-q}}{Z^V_n}\right) \notag \\
&  \leq \eta - \inf_{\lambda^\mathrm{ex} \in B_\delta(z)} \sum_{k=1}^{n-q} 
\mathcal J_V(\lambda_k) + \limsup_{n\to \infty}\,  (\beta'n)^{-1} \log \left(\frac{Z^{\frac{n}{n-q}V}_{n-q}}{Z^{V}_n}\right) .
\end{align}
Now we may apply Lemma \ref{lem:teknik} and use the fact that $\eta>0$ is arbitrary, to obtain
\begin{align} 
\limsup_{n\to \infty} \  (\beta'n)^{-1} \log \Pnv (\zeta_N^{(n)} \in B_\delta(z) ) \leq  -\inf_{\lambda^\mathrm{ex} \in B_\delta(z)} \sum_{k=1}^{q} 
\mathcal J_V(\lambda_k) + q \inf_{\xi\in \mathbb{R}} \mathcal J_V(\xi) . 
\end{align}
%Here we used that for $M$ large enough, $\mu_{V_M}=\mu_V$
Using that $\mathcal{J}_V$ is lower semicontinuous, the right hand side converges as $\delta\searrow 0$ to 
\begin{align}
- \sum_{k=1}^{q} \left( \mathcal J_V(z_k) - \inf_{\xi\in \mathbb{R}}\mathcal J_V(\xi)\right) 
=   - \mathcal{I}_N^\mathrm{ext}(z) . 
\end{align}
This concludes the proof of the upper bound.

\subsection{Proof of the lower bound}
\label{sect:extremalLDP3}

To prove the lower bound \eqref{lb}, we fix $z\in \mathcal{Z}_N$ and show that
\begin{align} \label{lb:lb}
\liminf_{n\to \infty}\ (\beta'n)^{-1} \log \Pnv(\zeta_N^{(n)} \in B_\delta(z)) \geq - \mathcal{I}_N^{\operatorname{ext}}(z) 
\end{align}
for $\delta$ small enough. 
We may restrict our proof to $z$ with $V(z_{i,j})<\infty$ for all $i,j$, as otherwise $\mathcal{I}_N^{\operatorname{ext}}(z)=\infty$ and the lower bound is trivial. We recall that in the proof of the upper bound, we identified $z$ with a vector in $\mathbb{R}^{q}$, $q=2MN$, and we decomposed the vector of $n$ eigenvalues as in \eqref{ub:representation0} into $\lambda^\mathrm{ex} \in \mathbb{R}^{q}$ and $\lambda^\mathrm{in} \in \mathbb{R}^{n-q}$. 
In the course of the proof, we will separate the $q$ extremal eigenvalues from the $n-q$ remaining eigenvalues and use the convergence of the empirical measure build from the latter ones. For this we need some care to separate the extremal eigenvalues from $\partial I$.  We will show \eqref{lb:lb} with $B_\delta(z)$ replaced by a set $U_\delta(z)\subset B_\delta(z)$, which is constructed as follows. $U_\delta(z)$ contains those $y\in \mathbb{R}^{q}$, where each coordinate $y_k$ deviates less than $\delta$ from $z_k$, if $z_k\notin \partial I$, and if $z_k\in \partial I$, $y_k$ keeps also a distance more than $\delta/2$ from $I$. More precisely, for $\delta>0$, let $J_\delta(x) = (x-\delta,x+\delta)$ if $x\notin I$, and if $x\in I$, let 
\begin{align} \label{lb:defB0}
J_\delta(x) = (x-\delta,x+\delta) \cap \{x': d(x',I)>\delta/2\} . 
\end{align}
Then, define
\begin{align} \label{lb:defB}
U_\delta(z) = J_\delta(z_1)\times \dots \times J_\delta(z_{q}). 
\end{align}
For $\delta$ small enough, $U_\delta(x)\subset I^c$, and additionally 
\begin{align*}
 \Pnv(\zeta_N^{(n)} \in B_\delta(z)) \geq  \Pnv(\zeta_N^{(n)}  \in U_\delta(z)). 
\end{align*}
Actually, to simplify later arguments, let $\delta$ be so small that $ x\in U_\delta(z)$ implies $d(x_k,I)>\delta/2$ for all $k$. Note that the latter condition is satisfied by definition for $k$ with $z_k\in \partial I$, but for the other only for $\delta$ small enough. Then we may bound
\begin{align}
\Pnv(\zeta_N^{(n)} \in U_\delta(z)) \geq \Pnv(\zeta_N^{(n)} \in U_\delta(z), d(\lambda_k)<\delta/4 \text{ for all } k> q) .
\end{align} 
Similar to \eqref{ub:representation}, we can then write
\begin{align} \label{lb:upsilon}
Z_n^V\, \binom{n}{q}^{-1} \Pnv(\zeta_N^{(n)} \in U_\delta(z), d(\lambda_k)<\delta/4 \text{ for all } k> q) = \int_{U_\delta(z)}  \Upsilon_{n,N}(\lambda^\mathrm{ex}) \diff \lambda^\mathrm{ex} , 
\end{align}
with $\Upsilon_{n,N}$ as in \eqref{ub:upsilon}, except that $\Xi_{n,N}(\lambda^\mathrm{ex})$ is replaced by 
\begin{align*}
\hat{\Xi}_{n,N}(\lambda^\mathrm{ex}) = \int_{D} \prod_{r=1}^{q}   \prod_{s=q+1}^{n} |\lambda_r-\lambda_s|^\beta  \prod_{r=q+1}^{n} e^{-n\beta'V(\lambda_r)}\prod_{q<  r <s \leq n} |\lambda_r - \lambda_s|^\beta  \diff \lambda^\mathrm{in} ,
\end{align*}
where the integration is now over the set $D=\{x :\, d(x_k,I)<\delta/4 \text{ for all } k\}$.
We then consider the probability measure $\chi_{n,N}$ on $\R^n$, which forces the extremal eigenvalues to be in $U_\delta(z)$, and is defined by
\begin{align} \label{lb:defchi}
\diff \chi_{n,N}(\lambda^\mathrm{ex}, \lambda^\mathrm{in}) := (\kappa_{n,N})^{-1} \left( \mathbbm{1}_{U_\delta(z)}(\lambda^\mathrm{ex})\diff \lambda^\mathrm{ex} \right) \left( \mathbbm{1}_{D}(\lambda^\mathrm{in})  \diff\mathbb P_{n-q}^{\frac{n}{n-q}V}(\lambda^\mathrm{in})\right) ,
\end{align}
where the arguments are $\lambda^\mathrm{ex}\in \mathbb{R}^{q}$ and $\lambda^\mathrm{in} \in \mathbb{R}^{n-q}$ %, $\lambda\mathrm{ex}$ is restricted to the set $K(x) = B(x)\cap \mathbb{R}^{\downarrow j} \cap \{V<\infty\}$, 
and $\kappa_{n,N}$ is the normalizing constant. Recall that as remarked at the end of Section \ref{sect:extremalLDP1}, we may assume that $V$ is infinite outside of a compact set.  
We then have the representation 
\begin{align} \label{lb:upsilon2}
\int_{U_\delta(z)}  \Upsilon_{n,N}(\lambda^\mathrm{ex}) \diff \lambda^\mathrm{ex} = Z_{n-q}^{\frac{n}{n-q}V} \kappa_{n,N} I_{n,N} ,
\end{align}
where, with $H(\lambda^\mathrm{ex})$ as in \eqref{ub:upsilon},
\begin{align}
I_{n,N} := \int 
H(\lambda^\mathrm{ex})  e^{-\beta'n\sum_{r=1}^{q} V(\lambda_r)}
 \left(\prod_{r=1}^{q}\prod_{s=q+1}^{n}|\lambda_r - \lambda_s|^\beta\right) \diff \chi_{n,N}(\lambda^\mathrm{ex}, \lambda^\mathrm{in}) . 
\end{align}
Jensen's inequality allows then to bound
\begin{align} \label{lb:Ibound}
\frac{1}{\beta'}\log I_{n,N}
\geq n I_{n,N}^{(1)} + 2I_{n,N}^{(2)} + 2 (n-q) I_{n,N}^{(3)} ,
\end{align}
where
\begin{align*}
I_{n,N}^{(1)} & = -\int \sum_{r=1}^{q} V(\lambda_r) \diff \chi_{n,N}(\lambda^\mathrm{ex}, \lambda^\mathrm{in}),\\
I_{n,N}^{(2)} & = \int \sum_{1\leq r<s \leq q} \log |\lambda_r-\lambda_s| \diff \chi_{n,N}(\lambda^\mathrm{ex}, \lambda^\mathrm{in}),\\
I_{n,N}^{(3)} & = \frac{1}{n-q}\int \sum_{r=1}^q \sum_{s=q+1}^{n} \log |\lambda_r - \lambda_s| \diff \chi_{n,N}(\lambda^\mathrm{ex}, \lambda^\mathrm{in}) .
\end{align*}
To obtain bounds for the $I_{n,N}^{(i)}$, we first consider the normalizing constant $\kappa_{n,N}$. From definition \eqref{lb:defchi}, it is given by
\begin{align*}
\kappa_{n,N} = \int_{U_\delta(z)}    \diff \lambda^\mathrm{ex} \, \mathbb{P}_{n-q}^{\frac{n}{n-q}V}(\lambda \in D).
\end{align*}
By definition of the set $D$, $I\subset \operatorname{Int}(D)$, and so by Lemma \ref{lem:cvproba}, 
\begin{align} \label{lb:kappa}
\lim_{n\to \infty} \mathbb{P}_{n-q}^{\frac{n}{n-q}V}(\lambda \in D) = 1,
\end{align}
which implies
\begin{align} \label{lb:kappa2}
\lim_{n\to \infty} \kappa_{n,N}  = |U_\delta(z)| , 
\end{align} 
where we write $|A|$ for the Lebesgue measure of a Borel set $A$. 
This allows then to prove the following limits for $I_{n,N}^{(i)}$, $i=1,2,3$:
\begin{align}
\lim_{n \to \infty} I_{n,N}^{(1)} & = -|U_\delta(z)|^{-1} \int_{U_\delta(z)} \sum_{r=1}^{q} V(\lambda_r) \diff \lambda^\mathrm{ex}, \\
\lim_{n \to \infty} I_{n,N}^{(2)} 
& =  |U_\delta(z)|^{-1} \int_{U_\delta(z)} \sum_{1\leq r<s \leq q} \log |\lambda_r-\lambda_s|  \diff \lambda^\mathrm{ex} \\
\lim_{n \to \infty} I_{n,N}^{(3)} & = 
|U_\delta(z)|^{-1} \int_{U_\delta(z)} \left( \sum_{r=1}^q \int \log |\lambda_r-\xi| \diff \mu_V(\xi) \right) \diff \lambda^{\mathrm{ex}}.
\end{align}
For the detailed arguments we again refer to \cite{magicrules}. This implies
\begin{align} \label{lb:Ibound2}
\liminf_{n\to \infty}  \frac{1}{\beta'n} \log I_{n,N} 
& \geq  -   |U_\delta(z)|^{-1} \int_{U_\delta(z)} \sum_{r=1}^q \mathcal{J}_V(z_r) \diff \lambda^\mathrm{ex} .
\end{align}
We can then return to \eqref{lb:upsilon} via \eqref{lb:upsilon2}, and obtain
\begin{align} \label{lb:almostfinished}
& \liminf_{n\to \infty}  \frac{1}{\beta'n} \log \Pnv(\zeta_N^{(n)} \in U_\delta(z), d(\lambda_k,I)<\delta/4 \text{ for all } k> q) \notag \\
& \geq -   |U_\delta(z)|^{-1} \int_{U_\delta(z)} \sum_{r=1}^q \mathcal{J}_V(\lambda_r) d\lambda^\mathrm{ex} 
		+ \liminf_{n\to \infty}  \frac{1}{\beta'n} \log \left( \frac{Z_{n-q}^{\frac{n}{n-q}V}}{Z_{n}^{V}} \right)
		+\liminf_{n\to \infty}  \frac{1}{\beta'n} \log \kappa_{n,N} .
\end{align}
By Lemma \ref{lem:teknik}, the first $\liminf$ in \eqref{lb:almostfinished} is given by $q \inf_{x \in \mathbb{R}} \mathcal{J}_V(x)$. Since $\kappa_{n,N}$ converges to a positive limit, the second one vanishes. 
Altogether, we obtain for $U_\delta(z)$ a neighborhood of a point $z\in \mathcal{Z}_N$ such that $\mathcal{I}_N^\ext(z)$ is finite, 
\begin{align} \label{lb:almostfinished2}
\liminf_{n\to \infty}  \frac{1}{\beta'n} \log \Pnv(\zeta_N^{(n)} \in B_\delta(x)) 
 \geq -  \left( |U_\delta(z)|^{-1} \int_{U_\delta(z)} \sum_{r=1}^q \mathcal{J}_V(\lambda_r) \diff \lambda^\mathrm{ex} - \inf_{x \in \mathbb{R}} \mathcal{J}_V(x) \right),
\end{align}
for $\delta >0$ small enough. Letting $\delta \to 0$, the set $U_\delta(z)$ concentrates at $z$ with $|U_\delta(z)|\to 0$. By continuity of $\mathcal{J}_V$ on the set where this function is finite, the lower bound in \eqref{lb:almostfinished2} converges to $\sum_k \mathcal{F}_V(z_k)$, which finishes the proof of the lower bound.

\appendix
\section{Appendix: Preliminaries on large deviations}

We consider a sequence $(X_n)_n$ of random variables with values in some Polish space $\mathcal{X}$ with Borel $\sigma$-algebra. Let $\mathcal{I}:\mathcal{X} \to [0,\infty]$ and $(a_n)_n$ be a sequence of positive real numbers with $a_n\to \infty$. We say that $(X_n)_n$ satisfies the large deviation principle with speed $a_n$ and rate function $\mathcal{I}$, if $\mathcal{I}$ is lower semicontinuous and 
\begin{itemize}
\item[(1)] for all $C\subset \mathcal{X}$ closed
\begin{align*}
\limsup_{n\to \infty} \frac{1}{a_n} \log \mathbb{P}(X_n \in C) \leq - \inf_{x\in C} \mathcal{I}(x) ,
\end{align*} 
\item[(2)] for all $O\subset \mathcal{X}$ open
\begin{align*}
\liminf_{n\to \infty} \frac{1}{a_n} \log \mathbb{P}(X_n \in O) \geq - \inf_{x\in C} \mathcal{I}(x) .
\end{align*} 
\end{itemize}
We will only consider good LDPs, which means that the level sets $\{\mathcal{I}\leq L\}$ are compact for all $L\geq 0$. The two following results are important tools that will be repeatedly applied in this paper. 

\medskip

\textbf{Contraction principle (Theorem 4.2.1 in \cite{demboz98}):} If $(X_n)_n$ satisfies the LDP in $\mathcal{X}$ with speed $a_n$ and good rate function $\mathcal{I}$ and if $f:\mathcal{X} \to \mathcal{Y}$ is a continuous mapping to another polish space $\mathcal{Y}$, then $(f(X_n))_n$ satisfies the LDP in $\mathcal{Y}$ with speed $a_n$ and good rate function given by
\begin{align*}
\mathcal{I}'(y) = \inf\{ \mathcal{I}(x) |\, f(x) = y\} .
\end{align*} 

\textbf{Dawson-G\"artner Theorem (Theorem 4.6.1 in \cite{demboz98}):} Let $X_n=(X_n^{(1)},X_n^{(2)},\dots )_n$ be a sequence of random variables with values in the sequence space $\mathcal{X}^\mathbb{N}$, equipped with the product topology, such that for any $j\geq 1$, the projection $(\pi_j(X_n))_n$ onto the first $j$ coordinates satisfies the LDP in $\mathcal{X}^j$ with speed $a_n$ and good rate function $\mathcal{I}_j$. Then $(X_n)_n$ satisfies the LDP with speed $a_n$ and good rate function 
\begin{align*}
\mathcal{I}(x) = \sup_{j\geq 1}\mathcal{I}(\pi_j(x)) .
\end{align*}

We also need the following technical result, proved in \cite{GNRaddendum}. It is tailored to the case when the rate is convex in one argument, but not in another one. Here, a function $F:\mathcal{X}\to (-\infty,\infty]$ is strictly convex at $x$, in which case $x$ is called an exposed point, if there exists $x^*$ in the topological dual $\mathcal{X}^*$ of $\mathcal{X}$, called an exposing hyperplane for $x$, such that 
\begin{align}\label{exposinghyper}
F(x) - \langle x^*,x\rangle < 
F(z) - \langle x^*,z\rangle
\end{align}
for all $z\neq x$.

\begin{thm}
\label{newgeneral}
Let $(X_n, Y_n)_n$ be an exponentially tight  sequence of random variables defined on $\mathcal X \times \mathcal Y$.
%Assume that $X_n \in \mathcal X$ and $Y_n \in \mathcal Y$ are defined on the same probability space and  
%that the two sequences $(X_n)_n$ and $(Y_n)_n$ are exponentially tight. 
Assume further that
\begin{enumerate}
\item There is a set $D\subset \mathcal X^*$ containing $0$ and functions $\Lambda:D\to \mathbb{R}$, $J:C_b(\mathcal Y) \to \mathbb{R}$ such that
for all $x^*\in D$ and $\varphi\in C_b(\mathcal Y)$
\begin{equation}
\label{newncgf}\lim_{n\to \infty} \frac{1}{n} \log \mathbb E \exp \left(n \langle x^* , X_n\rangle + n \varphi(Y_n)\right) =  \Lambda(x^*) + J (\varphi)\, .
\end{equation}
\item
If $\mathcal{F}$ denotes the set of exposed points $x$ of
\[\Lambda^* (x) = \sup_{x^* \in D} \{\langle x^*,x\rangle - \Lambda(x^*)\}   \]
with an exposing hyperplane $x^*$ satisfying $x^*\in D$ and 
$\gamma x^*\in D$ for some $\gamma% = \gamma(x)
>1 $, then for every $x\in \{ \Lambda^* < \infty\}$ there exists a sequence $(x_k)_k$ with $x_k \in \mathcal{F}$ such that $\lim_{k\to \infty} x_k = x$ and 
\begin{align*}
\lim_{k\to \infty} \Lambda^*(x_k) =\Lambda^*(x) .
\end{align*}
\end{enumerate}  
Then, the sequence $(X_n, Y_n)_n$ satisfies the LDP with speed $n$ and good rate function
\[\mathcal I(x,y) = \Lambda^*(x) + \mathcal{I}_Y(y)\,,\]
where 
\begin{equation*}
\mathcal{I}_Y(y) = \sup_{\varphi \in C_b(\mathcal Y)} \{ \varphi(y) - J(\varphi)\} .
\end{equation*}
\end{thm}

\bibliographystyle{alpha}
\bibliography{bibclean}

\end{document}